\documentclass[twocolumn,fleqn]{autart}    

\usepackage{graphicx}      
\usepackage{amsmath}
\usepackage{amssymb}
\usepackage{amsfonts}
\usepackage{macros_conf}
\usepackage{mathtools}
\usepackage{here}
\usepackage{macros_e}
\usepackage{mathfonts}
\usepackage{footnote}
\everymath{\displaystyle}
\renewcommand{\arraystretch}{0.95}
\DeclareMathOperator{\esssup}{ess \sup}
\def\GLTI{{\clG_\mathrm{LTI}}}
\def\GLTISI{{\clG_\mathrm{LTI,SI}}}
\def\wmax{{w_\mathrm{max}}}

\def\vabs{{v_\mathrm{abs}}}
\begin{document}
\begin{frontmatter}

\title{Induced Norm Analysis of Linear Systems for \\
Nonnegative Input Signals\thanksref{footnoteinfo}} 

\thanks[footnoteinfo]{%
A preliminary version of this paper was presented at ECC 2022 and IFAC WC 2023.  
This work was supported by JSPS KAKENHI Grant Number JP21H01354
and Japan Science and Technology Agency (JST) as part of Adopting Sustainable Partnerships for Innovative Research Ecosystem (ASPIRE), Grant Number JPMJAP2402.  
This work was also supported by the AI Interdisciplinary Institute ANITI funding, through the French "Investing for the Future PIA3" program under the Grant agreement n$^\circ$ANR-19-PI3A-0004 as well as by the National Research Foundation, Prime Minister's Office, Singapore under its Campus for Research Excellence and Technological Enterprise (CREATE) programme. 
}

\author[Kyushu]{Yoshio Ebihara} 
\author[Kyukou]{Noboru Sebe} 
\author[IMI]{Hayato Waki} 
\author[LAAS]{Dimitri Peaucelle}
\author[LAAS]{Sophie Tarbouriech}
\author[LAAS]{Victor Magron}
\author[Kyoto]{Tomomichi Hagiwara} 

\address[Kyushu]{Faculty of Information Science and Electrical Engineering, 
Kyushu University, 
Fukuoka 819-0395, Japan.  \\
(e-mail: ebihara@ees.kyushu-u.ac.jp).}
\address[Kyukou]{Department of Intelligent and Control Systems,
Kyushu Institute of Technology, Fukuoka 820-8502, Japan.  }
\address[IMI]{Institute of Mathematics for Industry,
Kyushu University, 
Fukuoka 819-0395, Japan.  }
\address[LAAS]{LAAS-CNRS, Universit\'{e} de Toulouse, CNRS, F-31400, Toulouse, France.  }
\address[Kyoto]{Department of Electrical Engineering,
Kyoto University, 
Kyoto 615-8510, Japan.  
}


\begin{abstract}                
This paper is concerned with the analysis of the $L_p\ (p\in[1,\infty),\ p=\infty)$ induced norms 
of continuous-time linear systems where input signals are restricted to be nonnegative.  
This norm is referred to as the $L_{p+}$ induced norm in this paper.  
It has been shown recently that the $L_{2+}$ induced norm is effective for 
the stability analysis of nonlinear feedback systems where the nonlinearity
returns only nonnegative signals.  
However, the exact computation of the $L_{2+}$ induced norm is essentially difficult.  
To get around this difficulty, in the first part of this paper, we provide a 
copositive-programming-based method for the upper bound computation 
by capturing the nonnegativity of the input signals by copositive multipliers.  
In the second part, 
we consider how far the $L_{2+}$ induced norm can be smaller 
than the standard $L_{2}$ induced norm, and
derive the uniform infimum on the ratio of
the $L_{2+}$ induced norm to the $L_{2}$ induced norm 
over all linear systems including infinite-dimensional ones.  
Then, for each linear system, we finally derive a computation method of the 
lower bounds of the $L_{2+}$ induced norm 
that are larger than (or equal to) the value determined by the uniform infimum.  
The effectiveness of the upper/lower bound computation methods is illustrated by numerical examples.  
\end{abstract}

\end{frontmatter}

\section{Introduction}

The $L_{p}\ (p\in[1,\infty),\ p=\infty)$ induced norm plays an important role
in dynamical system analysis \cite{Desoer_1975,Khalil_2002}.  
The induced norm is the core in assessing
the input-output stability ($L_p$ stability) of dynamical systems, 
and also serves as a reasonable measure for disturbance attenuation.  
The induced norm is also particularly useful in investigating the stability 
of interconnected systems via the small-gain theorem \cite{Khalil_2002}.  
In this study, we introduce the $L_p\ (p\in[1,\infty),\ p=\infty)$ induced norms 
for continuous-time linear systems where input signals are restricted to be nonnegative.  
This norm is referred to as the $L_{p+}$ induced norm in this paper.  

Recently, there has been a growing attention on 
control theoretic approaches for the analysis and synthesis 
of optimization algorithms \cite{Lessard_SIAM2016}, 
static feedforward neural networks (NNs) 
\cite{Raghunathan_NIPS2018,Raghunathan_ICLR2018,Fazlyab_arxiv2019,Fazlyab_IEEE2022,Chen_NIPS2020,Ebihara_ECC2024}, 
dynamical NNs such as recurrent NNs (RNNs) 
\cite{Revay_LCSS2021,Scherer_IEEEMag2022},  
and dynamical systems driven by NN controllers \cite{Yin_IEEE2022}.   
By capturing the input-output behavior of nonlinearities 
in the algorithms or NNs via quadratic constraints,
we can cast the analysis and synthesis problems    
into numerically tractable semidefinite programming problems (SDPs).  
Along this stream, in \cite{Ebihara_EJC2021,Ebihara_CDC2021,Motooka_ISCIE2022}, 
we dealt with the stability analysis of RNNs 
with activation functions being rectified linear units (ReLUs).  
In particular, by focusing on the fact that the ReLUs return only nonnegative signals, 
we derived an $L_{2+}$-induced-norm-based small-gain theorem \cite{Ebihara_CDC2021,Motooka_ISCIE2022}.  
Still, exact (or smaller upper bound) computation of the $L_{2+}$ induced norm 
of linear systems 
remained to be an outstanding issue.  

This study is also motivated by recent advancement on the study of positive systems
(\cite{Blanchini_IEEE2012,Briat_IJRN2013,Tanaka_IEEE2011,Blanchini_2015,Rantzer_EJC2015,Rantzer_IEEE2016,Rantzer_ARC2021,Ebihara_IEEE2017,Kato_LCSS2020}).   
A positive system is a system whose 
state variables and/or output variables are nonnegative.  
This property arises
naturally in biology, network communications, economics, and
probabilistic systems, and representative 
linear time-invariant (LTI)
positive systems include compartmental systems. Moreover, simple LTI systems
such as integrators and first-order lags and their series/parallel
connections are all positive, and these are often employed as
typical models of moving objects. Even though their dynamics
are very simple, the behavior of large-scale multi-agent systems constructed
from them is complicated and deserves investigation.  
Still, the system positivity allows us to employ Lyapunov functions that are 
linear or quadratic and diagonal with respect to the state,  
and this enables us to derive analysis and synthesis conditions 
that scale linearly to the system size \cite{Briat_IJRN2013,Tanaka_IEEE2011,Rantzer_EJC2015,Rantzer_IEEE2016,Ebihara_IEEE2017}.  
On the other hand, in \cite{Kato_LCSS2020}, we recognized the usefulness
of the copositive multipliers and the copositive programming problems (COPs) \cite{Dur_2010}
in capturing the signal nonnegativity in standard quadratic fashion.  
We thus had a prospect that the COPs can be used
even for nonpositive (general) linear system analysis 
in handling nonnegative signals.  

On the basis of the preceding studies, 
in this paper, we investigate 
the $L_{p+}\ (p\in[1,\infty),\ p=\infty)$ induced norms 
of continuous-time linear systems from a broad perspective.  
Since the upper bound of 
$L_{2+}$ induced norm is especially important in 
performance guarantees of dynamical systems,  
we first consider its upper bound computation problem (Problem~1).  
By introducing positive filters of increased degree and
employing copositive multipliers of increased size (freedom), 
we derive a sequence of COPs that generates a
monotonically non-increasing sequence of upper bounds.   
Furthermore, by applying an inner approximation to the copositive cone,
we derive a numerically tractable sequence of SDPs.  
Then, we second tackle the 
analysis of the uniform infimum on the ratio of
the $L_{2+}$ induced norm to the $L_{2}$ induced norm 
over all linear systems including infinite-dimensional ones (Problem~2).  
By definition, $L_{p+}$ induced norm is 
smaller than (or equal to) the standard $L_p$ induced norm,   
and Problem~2 is motivated to
clarify how far the $L_{p+}$ induced norm can be smaller than the 
$L_p$ induced norm.  
Here we deal with $L_{p+}$ induced norms for $p\in[1,\infty)$ and $p=\infty$ 
in a unified fashion, 
since the underlying methodology is independent of $p$.  
The result of Problem~2 is applicable to any LTI systems 
to obtain their lower bounds of the $L_{2+}$ induced norm.  
However, for each linear system, 
it is expected that we can obtain
better (larger) lower bounds than the value determined by the uniform infimum.  
Such a lower bound for the $L_{2+}$ induced norm 
is desirable to evaluate the 
accuracy of the upper bounds obtained for Problem~1.  
We therefore finally 
deal with the lower bound analysis problem of the $L_{2+}$ induced norm for a given
linear system (Problem~3).  
By reducing the lower bound analysis problem into 
a semi-infinite programming problem  \cite{Shapiro_2009},     
we derive a computation method that enables us to obtain
lower bounds that are larger than (or equal to) the value determined by the uniform infimum.  

A preliminary version of this paper is published in \cite{Ebihara_ECC2022,Ebihara_IFAC2023}, 
where the above three problems have been investigated.  
The novel aspects of the present paper over  \cite{Ebihara_ECC2022,Ebihara_IFAC2023}
are thoroughly summarized in \rsub{sub:novel}, 
after more accurate problem statements in \rsub{sub:problem}.  

\section{Preliminaries, Motivations, and Problem Settings}

\subsection{Notation}

We use the following notation in this paper.  
The set of natural numbers is denoted by $\bbN$.  
The set of $n$-dimensional real vectors 
(with nonnegative entries) is denoted by $\bbR^n\ (\bbR_+^n)$, 
and the set of $n\times m$ real matrices  (with nonnegative entries) 
is denoted by $\bbR^{n\times m}\ (\bbR_+^{n\times m})$.  
For $M\in\bbR^{n\times m}$, we also write $M\ge 0$ 
to denote $M\in \bbR_+^{n\times m}$.   
The set of $n\times n$ real symmetric (positive definite) 
matrices is denoted by $\bbS^n\ (\bbS_{++}^n)$.  
The set of $n\times n$ Hurwitz and Metzler matrices are denoted by 
$\bbH^n$ and $\bbM^n$, respectively, 
where a matrix $M\in\bbR^{n\times n}$  is said to be Metzler
if $M_{ij}\ge 0\ (i\ne j)$.  
For $S\in\bbS^n$, we write $S\succeq 0\ (S\preceq 0)$ to
denote that $S$ is positive (negative) semidefinite.  
For $w_1,w_2\in \bbR^n$, we define
$\wmax=\max(w_1,w_2)\in \bbR^n$ by 
$w_{\mathrm{max},i}=\max(w_{1,i},w_{2,i})\ (i=1,\cdots,n)$.  
%

\subsection{Preliminaries on Signals, Norms, and Cones}
For a vector $v\in\bbC^{n_v}$, we define its $l_p\ (p\in[1,\infty),\ p=\infty)$ norm by
\[
\begin{aligned}
& |v|_p:=\left(\sum_{i=1}^{\nv}  |v_i|^p\right)^{1/p}\!\!\! (p\in[1,\infty)),\  |v|_\infty:=\!\!\!\max_{i=1,\cdots,\nv} |v_i|.  
\end{aligned}
\]
For a matrix $M\in \bbC^{n\times m}$, we define its $l_p\ (p\in[1,\infty),\ p=\infty)$ induced norm by
\[
 \|M\|_p=\max_{v\in\bbC^m,\ |v|_p=1} |Mv|_p\ (p\in[1,\infty),\ p=\infty).  
\]
For a continuous-time signal $w:\ [0,\infty) \to \bbR^{\nw}$, we 
define its $L_p\ (p\in[1,\infty),\ p=\infty)$ norm by
\[
\|w\|_{p}:=\left(\int_{0}^{\infty}|w(t)|_p^p dt\right)^{1/p}\ (p\in[1,\infty)),
\]
\[
\|w\|_\infty:=\esssup\displaylimits_{0\le t<\infty}|w(t)|_\infty. 
\]
For $p\in [1,\infty)$ and $p=\infty$, we also define
the (standard) $L_p$ space and the set of the nonnegative $L_p$ signals by
\[
\begin{array}{@{}l}
 L_{p}  :=\left\{w:\ \|w\|_{p}<\infty \right\},\\ 
 L_{p+} :=\left\{w:\ w\in L_{p},\ w(t)\ge 0\ \mbox{\rm {a.e. in}}\ t\in[0,\infty)\right\}.    
\end{array}
\]
We also define their extended versions
\[
\begin{array}{@{}l}
 L_{p,\mathrm{e}}  :=\left\{w:\ w_\tau\in L_{p}\ \forall \tau\in[0,\infty) \right\},\\ 
 L_{p,\mathrm{e}+} :=\left\{w:\ w\in L_{p,\mathrm{e}},\ w(t)\ge 0\ \mbox{\rm {a.e. in}}\ t\in[0,\infty)\right\}.    
\end{array}
\]
where $w_\tau$ is a truncation of $w:\ [0,\infty) \to \bbR^{\nw}$ defined by
\[
 w_\tau(t) = \left\{
\begin{array}{cc}
 w(t)& 0 \le t \le \tau, \\
 0    & \tau<t. \\
\end{array}
\right.
\]
For a linear operator 
\begin{equation}
G:\ L_{p}\ni w \mapsto z \in L_{p}\ (p\in[1,\infty),\ p=\infty),  
\label{eq:G}
\end{equation}
we define its (standard) $L_p$ induced norm by
\[
\|G\|_{p}:=\sup_{w\in L_p,\ \|w\|_p=1} \ \|z\|_p.  
\]
On the other hand, we newly introduce 
\begin{equation}
\|G\|_{p+}:=\sup_{w\in L_{p+},\ \|w\|_p=1} \ \|z\|_p.  
\label{eq:Lp+ind}
\end{equation}
%
This obviously satisfies the axioms of a norm and is 
referred to as the $L_{p+}$ induced norm in this paper.  
Due to the restriction to nonnegative input signals, 
we can readily see that the very basic property $\|G\|_{p+}\le \|G\|_{p}\ (p\in[1,\infty),\ p=\infty)$ holds.  

We define the positive semidefinite cone $\PSD^{n}\subset \bbS^n$, 
the copositive cone $\COP^{n}\subset \bbS^n$, and the nonnegative cone
$\NN^{n}\subset \bbS^n$ as follows:
\[
 \begin{array}{@{}l}
  \PSD^n:=\{P\in\bbS^n:\ x^TPx\geq 0\ \forall x\in\bbR^n\},\\
  \COP^n:=\{P\in\bbS^n:\ x^TPx\geq 0\ \forall x\in\bbR_{+}^n\},\\
  \NN^n:= \{P\in\bbS^n:\ P\geq 0\}.  
 \end{array}
\]
The semidefinite programming problem (SDP) and 
the copositive programming problem (COP) 
are convex optimization problems in which 
we minimize a linear objective function over the 
linear matrix inequality (LMI) constraints 
on $\PSD$ and $\COP$, respectively.     
As mentioned in \cite{Dur_2010}, the COP is 
a co-NP complete problem and hence numerically intractable in general.  
However, the convex optimization problems on $\PSD+\NN$ 
with ``$+$'' being the Minkowski sum is essentially an SDP and hence numerically tractable.  
Since $\PSD\subset \PSD+\NN\subset \COP$ obviously holds, 
we can therefore apply an inner approximation to $\COP$ with 
$\PSD+\NN$ and solve a COP in a sufficient fashion.  
In particular, since 
$\PSD^n+\NN^n=\COP^n\ (n\le 4)$ holds \cite{Dur_2010},
we can solve a COP exactly if $n\le 4$.  

\subsection{Motivations: Relevance of $L_{p+}$ Induced Norm to Dynamical System Analysis and Synthesis}
\label{sub:motiv}

This subsection provides two typical examples that motivate us
to focus on $L_{p+}$ induced norms of dynamical systems.  

\subsubsection{Evaluation of Difference of Positive Systems}
\label{sub:diff_pos}

We first recall the definition of positive systems and related results.
\begin{definition}[\cite{Farina_2000,Kaczorek_2001}]
Let us consider the continuous-time finite-dimensional LTI system $G$ given by
\begin{equation}
G:\ 
\left\{
\arraycolsep=0.5mm
\begin{array}{cccccccc}
 \dot x(t)&=& A x(t)& + &B w(t), \\
   z(t)&=& C x(t)& + &D w(t) \\
\end{array}
\right.  
\label{eq:GLTI}
\end{equation}
where
$A\in\bbR^{n\times n}$, 
$B\in\bbR^{n\times \nw}$,
$C\in\bbR^{\nz\times n}$, and
$D\in\bbR^{\nz\times \nw}$.   
Then, the system $G$ given by \rec{eq:GLTI} is said
to be {\it externally positive} if its output is nonnegative for any nonnegative
input under zero initial state.  
In addition, it is said
to be {\it internally positive} if its state and output are nonnegative for 
any nonnegative input and nonnegative initial state.  
\end{definition}
\begin{proposition}[\cite{Farina_2000,Kaczorek_2001}]
The system $G$ given by \rec{eq:GLTI}
is externally positive if and only if its impulse response is nonnegative.   
In addition, it is internally positive if and only if 
$A\in\bbM^{n}$, $B\in\bbR_+^{n\times \nw}$, 
$C\in\bbR_+^{\nz\times \nw}$, and $D\in\bbR_+^{\nz\times \nw}$.  
\label{pr:positive}
\end{proposition}

By definition, it is true that if $G$ is internally positive 
then it is externally positive.  
It is also well known that if $G$ is externally positive then 
$\|G\|_2=\|G\|_{2+}$ holds \cite{Tanaka_IEEE2011,Rantzer_IEEE2016}.    
This result can readily be generalized in the following way.  
\begin{proposition}
Suppose $G$ is externally positive. 
Then we have $\|G\|_p=\|G\|_{p+}\ (p\in[1,\infty),\ p=\infty)$.   
\label{pr:expos_result}
\end{proposition}
The proof of this proposition is given in Appendix~\ref{sec:ap0}
to make the paper self-contained.    

On the basis of the above preliminaries, 
let us now consider two externally positive systems $G_1$ and $G_2$
of the same input-output size.  
Here we want to evaluate the ``difference'' (or say, error) between them.  
This issue typically arises in 
positivity-preserving model reduction for positive systems \cite{Reis_2009,Li_Automatica2011,Sootla_CDC2012}.  
A common and easy-to-compute 
measure is the $L_2$ induced norm
($H_\infty$ norm) of the error system, i.e., $\|G_1-G_2\|_2$.  
However, since the input of positive systems are often naturally 
nonnegative, it is more suitable to evaluate the difference
under nonnegative inputs.  This leads to the requirement to evaluate
$\|G_1-G_2\|_{2+}$ (or $\|G_1-G_2\|_{p+}$ for broader treatments).  
It should be noted that the error system $G_1-G_2$ is no longer 
externally positive in general even if 
$G_1$ and $G_2$ are.  
Therefore the computation of $\|G_1-G_2\|_{2+}$ is
not trivial and becomes a challenging issue.  

\subsubsection{Stability Analysis of Dynamical Systems Driven by Neural Network Controllers}
\label{sub:ReLU-Feedback}

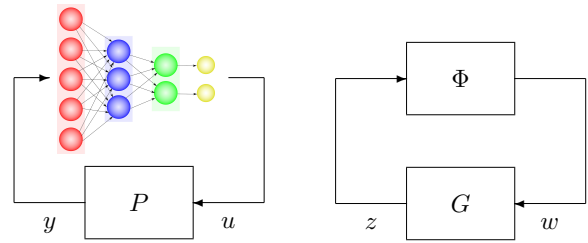
\begin{figure}[t]
\begin{center}
\begin{center}
\scalebox{0.95}{$
\hspace*{3mm}
\begin{minipage}{45mm}
\hspace*{5mm}
\includegraphics[width=2.2cm]{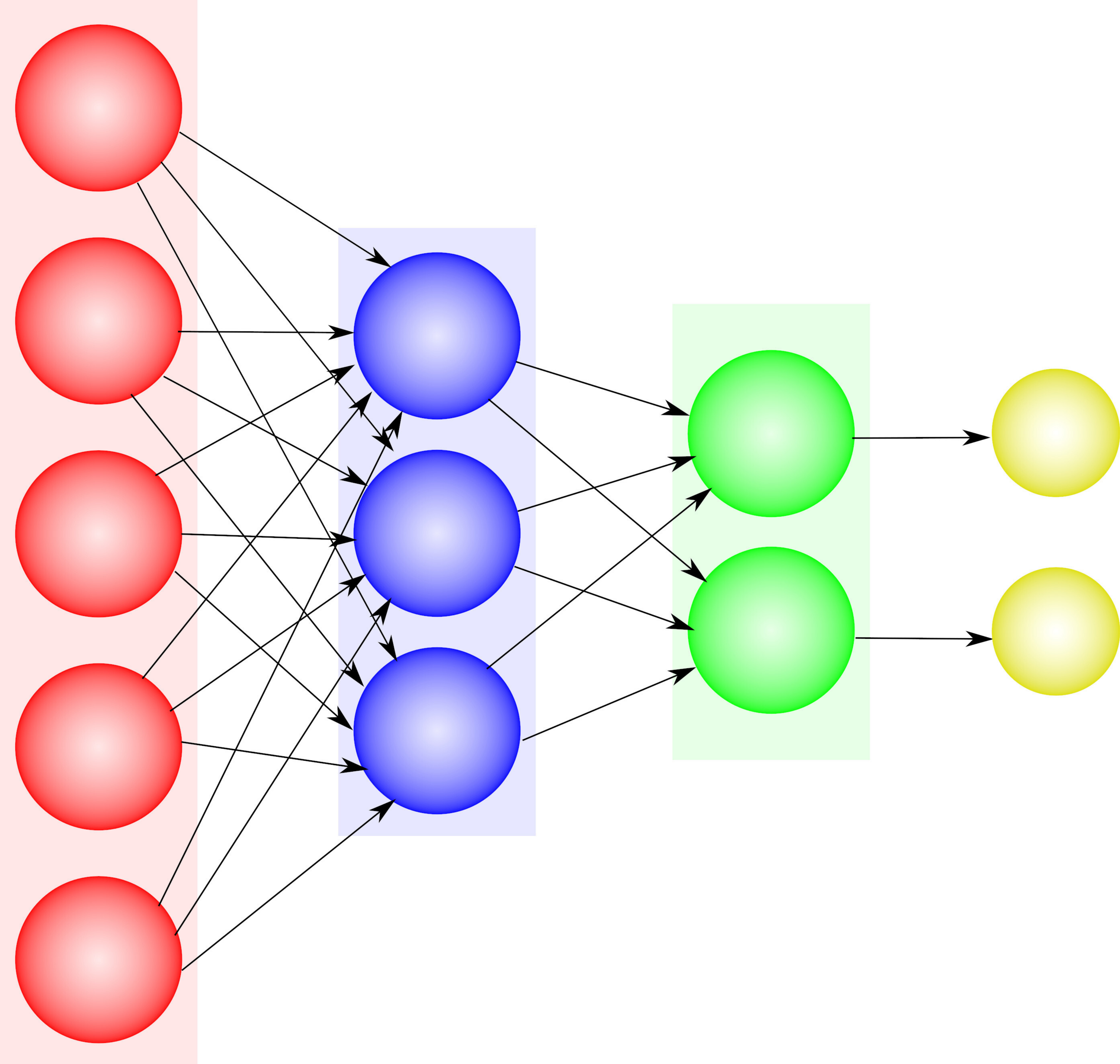}\vspace*{-11.0mm}\\
\begin{picture}(3.5,2.25)(0,0)
\put(0,2.25){\vector(1,0){0.5}}
\put(3,2.25){\line(1,0){0.5}}
\put(3.5,2.25){\line(0,-1){1.75}}
\put(3.5,0.5){\vector(-1,0){1.0}}
\put(3,0.3){\makebox(0,0)[t]{$u$}}
\put(1,0){\framebox(1.5,1){$P$}}
\put(1,0.5){\line(-1,0){1.0}}
\put(0.5,0.3){\makebox(0,0)[t]{$y$}}
\put(0,0.5){\line(0,1){1.75}}

\end{picture} 
\end{minipage}
\begin{minipage}{40mm}
\vspace*{5.5mm}
\begin{picture}(3.5,2.75)(0,0)
\put(0,2.25){\vector(1,0){1.0}}
\put(1,1.75){\framebox(1.5,1){$\Phi$}}
\put(2.5,2.25){\line(1,0){1.0}}
\put(3.5,2.25){\line(0,-1){1.75}}
\put(3.5,0.5){\vector(-1,0){1.0}}
\put(3,0.3){\makebox(0,0)[t]{$w$}}
\put(1,0){\framebox(1.5,1){$G$}}
\put(1,0.5){\line(-1,0){1.0}}
\put(0.5,0.3){\makebox(0,0)[t]{$z$}}
\put(0,0.5){\line(0,1){1.75}}

\end{picture} 
\end{minipage}$}
\end{center}
\caption{Dynamical Systems Driven by NN Controllers (left) and Nonlinear Feedback System Representation (right).  }
\label{fig:NNfb}
\end{center}
\vspace*{-1mm}
\end{figure}

Recently, control theoretic approaches have attracted great attention  
for the analysis of static feedforward neural networks (NNs) 
\cite{Raghunathan_NIPS2018,Raghunathan_ICLR2018,Fazlyab_arxiv2019,Fazlyab_IEEE2022,Chen_NIPS2020,Ebihara_ECC2024}, 
dynamical NNs such as recurrent NNs 
\cite{Scherer_IEEEMag2022,Revay_LCSS2021,Ebihara_EJC2021,Ebihara_CDC2021,Motooka_ISCIE2022}, 
and dynamical systems driven by 
NN controllers \cite{Yin_IEEE2022} as shown in \rfig{fig:NNfb} (left).  
Suppose there are no external signals 
(such as bias signals) to 
the nonlinear activation functions in the NNs.  
Then, the dynamical NNs and linear dynamical systems driven by NNs can be modeled 
as a nonlinear feedback system shown in \rfig{fig:NNfb}  (right), 
where $G$ is typically a stable LTI system and 
$\Phi:\ \bbR^m\to \bbR^m$ is a static nonlinear operator 
representing nonlinear activation functions in the NNs.  
Such modeling is useful if we are interested in
the global asymptotic stability of the origin (GAS) of the original systems.  
In particular, in the case where all the activation functions employed are 
rectified linear units (ReLUs), we see that
$\Phi:\ \bbR^m\to \bbR^m$ satisfies 
$\|\Phi\|_2=1$ and $\Phi:\ \bbR^m\to \bbR_+^m$, i.e., 
$\Phi$ returns only nonnegative signals.  
In such a typical case, the $L_{2+}$ induced norm becomes quite relevant
for the stability analysis of the nonlinear feedback system as briefly explicated below.  

From the standard $L_2$-induced-norm-based 
small-gain theorem \cite{Khalil_2002}, 
we see that the feedback system shown in \rfig{fig:NNfb} is
(well-posed and) GAS if $\|G\|_2<1$.  
On the other hand, by actively using the nonnegativity nature of $\Phi$, 
it has been shown recently in \cite{Motooka_ISCIE2022} that 
the feedback system is (well-posed and) GAS if $\|G\|_{2+}<1$.  
As illustrated by this concrete example, 
the $L_{2+}$-induced-norm-based small-gain theorem
has potential abilities for the stability analysis of feedback systems
with nonnegative nonlinearities.   
However, exact computation of the $L_{2+}$ induced norm 
is inherently difficult and this strongly motivates the current study
on its upper and lower bounds computation.   
The signal-nonnegativity-based analysis of NNs 
can also be found at \cite{Gronqvist_MTNS2022}.  

\subsection{Problem Settings}
\label{sub:problem}

This subsection poses three problems to be investigated in this paper.  
Since the upper bound of $L_{2+}$ induced norm is important in guaranteeing the system performance
as illustrated in \rsub{sub:motiv}, 
we naturally consider the next problem.  
\begin{problem}[$L_{2+}$ Upper Bound Analysis]\ \\
For a given stable LTI system $G$ of the form \rec{eq:G}, 
find an upper bound of $\|G\|_{2+}$ that is as small as possible.  
\end{problem}

On the other hand, $\|G\|_{p+}\le \|G\|_p$ does hold for 
$p\in[1,\infty)$ and $p=\infty$.  
On this relationship, 
we are interested in
how far $\|G\|_{p+}$ can be smaller than $\|G\|_p$
partly from an applied mathematics point of view.  
To clarify this point, we investigate the next problem.  
\begin{problem}[Uniform Infimum of $\|G\|_{p+}/\|G\|_p$]
For each $p\in[1,\infty)$ and $p=\infty$, find the uniform infimum on the ratio of
the $L_{2+}$ induced norm to the $L_{2}$ induced norm defined by
$\nu_p^\star:=\inf_{G\in \GLTI} \|G\|_{p+}/\|G\|_p$,    
which is alternatively characterized as
\begin{equation}
 \nu_p^\star=\sup\left\{\nu_p\in\bbR:\ \|G\|_{p+}\ge \nu_p \|G\|_{p}\ \forall G\in \GLTI\right\}.  \hspace*{-10mm}
\label{eq:nups} 
\end{equation}
Here, $\GLTI$ stands for the set of stable and causal LTI systems including infinite dimensional ones.  
\end{problem}
\begin{remark}
We use the terminology ``uniform'' since $\nu_p^\star$ is valid for any $G\in\GLTI$.  
From \rec{eq:nups}, it is clear that 
$\nu_p\|G\|_p$ is a lower bound of $\|G\|_{p+}\ (p\in[1,\infty),\ p=\infty)$ for any $G\in\GLTI$.  
\end{remark}

The result of Problem~2 is applicable to any $G\in\GLTI$ to obtain a lower bound of $\|G\|_{2+}$.  
However, for each $G\in\GLTI$, it is expected that we can obtain
better (larger) lower bounds than the value determined by the uniform infimum.  
Such a lower bound for the $L_{2+}$ induced norm is desirable to evaluate the 
accuracy of the upper bound obtained for Problem~1.  Therefore we consider the next problem.  
\begin{problem}[$L_{2+}$ Lower Bound Analysis]\ \\
For a given stable LTI system $G$ of the form \rec{eq:G}, 
find a lower bound of $\|G\|_{2+}$ that is as large as possible.  
\end{problem}
%

\subsection{Novel Aspects Over \cite{Ebihara_ECC2022,Ebihara_IFAC2023}}
\label{sub:novel}

The novel aspects of the present paper over our preceding studies in \cite{Ebihara_ECC2022,Ebihara_IFAC2023} 
are summarized as follows:
\begin{quote}
\begin{itemize}
 \item[(i)] Problem~1 is dealt with in \cite{Ebihara_ECC2022}.  
	    However, one of the key results in \cite{Ebihara_ECC2022}, i.e., 
	    the construction of the monotonically non-increasing sequence of the upper bounds
	    (\rth{th:monotone} in this paper), is given without any proof.  
	    The explicit proof is now given in \rsec{sec:ub}, thereby the proposed 
	    upper bound computation method is completed for the first time in this paper.   
 \item[(ii)] Problem~3 is investigated in \cite{Ebihara_IFAC2023}.  
	     It is nonetheless true that the results there are restricted 
	     to be single-input systems due to technical difficulties.  
	     In the present paper we succeed in extending the results to multi-input systems (\rth{th:lbcomp_multi}).  
\item[(iii)] The effectiveness of the upper bound computation method for Problem~1
	     and the lower bound computation method for Problem~3 is illustrated
	     for multi-input systems in \rsub{sub:num}.  
	     In particular, we provide an interesting numerical example on the evaluation of
	     the difference of positive systems.  
\end{itemize}
\end{quote}
%


\section{Main Results for $L_{2+}$ Upper Bound Analysis Problem}
\label{sec:ub}

\subsection{Upper Bound Computation of $L_{2+}$ Induced Norm by Positive Filters}
\label{sub:ub}

Let us focus on the case where $G$ in \rec{eq:G}
is a finite-dimensional LTI system given by \rec{eq:GLTI} where $x(0)=0$.  
We assume that the system $G$ is stable, i.e., $A\in\bbH^n$.  
As noted previously, it is very clear that $\|G\|_{2+}\le \|G\|_{2}$, 
i.e., $\|G\|_{2}$ is a trivial upper bound of $\|G\|_{2+}$.  

For better upper bound computation of $\|G\|_{2+}$, 
it is promising to actively use the fact that the input signal
$w$ is restricted to be nonnegative.    
To this end, let us introduce the positive filter given by
\begin{equation}
G_p : \begin{cases}
 \dot x_p(t)=A_p x_p(t)+B_pw(t),\ x_p(0)=0,\\
 z_p(t)= 
 \left[
  \begin{array}{c}
   I_{\np}  \\
   0_{\nw,\np}
  \end{array}
 \right] x_p(t)+
 \left[
  \begin{array}{c}
   0_{\np,\nw}\\
   I_{\nw}  \\
  \end{array}
 \right] w(t)\hspace*{-20mm}
\end{cases}
\label{eq:Gp}
\end{equation}
where $A_p\in \bbH^{\np}\cap \bbM^{\np}$, $B_p\in\bbR_+^{\np\times \nw}$.  
It is clear that the filter $G_p$ is internally positive from \rpr{pr:positive}.  

We next consider the vertical stacking of $G$ above $G_p$ 
and construct the augmented system $G_a$ given by
\begin{equation}
 \scalebox{0.92}{$
 \begin{array}{@{}l}
G_a : \begin{cases}
\dot x_a(t) = A_a x_a(t)+ B_a w(t),\\
z(t)= C_a x_a(t) +D_aw(t),\\
 z_p(t)= 
 \left[
  \begin{array}{cc}
   0_{\np,n} & I_{\np}  \\
   0_{\nw,n} & 0_{\nw,\np}
  \end{array}
 \right]x_a(t)+
 \left[
  \begin{array}{c}
   0_{\np,\nw}\\
   I_{\nw}  \\
  \end{array}
 \right] w(t),  \hspace*{-20mm}
\end{cases}\\
x_a:=
\left[
  \begin{array}{c}
   x \\
   x_p \\
  \end{array}
\right],\ 
A_a:= 
\left[
  \begin{array}{cc}
   A & 0  \\
   0 & A_p \\
  \end{array}
\right],\ 
B_a:=
\left[
  \begin{array}{c}
   B \\
   B_p \\
  \end{array}
\right],\\
C_a:=
\left[
  \begin{array}{cc}
   C & 0_{n_z,\np}  \\
  \end{array}
\right], D_a:=D.  
\end{array}$}
\label{eq:Ga}
\end{equation}
In the above augmented system $G_a$, 
it is very important to note that 
the output $z_p$ is nonnegative for any nonnegative 
input $w$.  
By focusing on this nonnegativity property, 
the next result can be obtained.  
The proof of this theorem is given in 
Appendix~\ref{sec:ap1}.  
\begin{theorem}
For the LTI system $G$ given by \rec{eq:GLTI} and a given $\gamma>0$, 
we have $\|G\|_{2+}\le \gamma$ if 
there exist $P_a \in \bbS^{n+\np}$ and $Q_a\in \COP^{\np+\nw}$ such that
\begin{equation}
\scalebox{0.9}{$
\begin{array}{@{}l}
\begin{bmatrix}
P_aA_a + A_a^T P_a + C_a^T C_a & P_aB_a + C_a^T D_a \\
B_a^T P_a + D_a^T C_a & D_a^T D_a- \gamma ^2 I_{\nw}
\end{bmatrix}\\
+
\begin{bmatrix}
0_{n,\np+\nw} \\ I_{\np+\nw}
\end{bmatrix}
Q_a
\begin{bmatrix}
0_{n,\np+\nw} \\ I_{\np+\nw}
\end{bmatrix}^T
\preceq 0.  
\end{array}$}
\label{eq:L2+COPnew}
\end{equation}
\label{th:main1}
\end{theorem}
We note that if we do not employ \rec{eq:Gp} 
(i.e., if we consider the \textit{filter-free} case), 
the condition \rec{eq:L2+COPnew} reduces to
\begin{equation}
\scalebox{0.9}{$
\begin{array}{@{}l}
\begin{bmatrix}
PA + A^T P + C^T C & PB + C^T D \\
B^T P + D^T C & D^T D- \gamma ^2 I_{\nw}
\end{bmatrix}\\
+
\begin{bmatrix}
0_{n,\nw} \\ I_{\nw}
\end{bmatrix}
Q
\begin{bmatrix}
0_{n,\nw} \\ I_{\nw}
\end{bmatrix}^T
\preceq 0,\ 
P \in \bbS^{n},\ Q\in \COP^{\nw}.  
\end{array}$}
\label{eq:L2+COPnewFilterFree}
\end{equation}
On the basis of \rth{th:main1}, let us consider the following 
COP and SDP:  
\begin{equation}
 \scalebox{1.0}{$
\begin{array}{@{}l}\displaystyle
\ogam_a := \inf_{\gamma, P_a, Q_a}\ \gamma\quad \mathrm{subject\ to}\  
\rec{eq:L2+COPnew},\\ 
P_a \in \bbS^{n+\np},\ Q_a \in \COP^{\np+\nw}, 
\end{array}$}
\label{eq:COPnew}
\end{equation}
\begin{equation}
 \scalebox{1.0}{$
\begin{array}{@{}l}\displaystyle
\oogam_a := \inf_{\gamma, P_a, Q_a}\ \gamma\quad \mathrm{subject\ to}\
 \rec{eq:L2+COPnew},\\ 
P_a \in \bbS^{n+\np},\ Q_a \in \PSD^{\np+\nw}+\NN^{\np+\nw}.  
\end{array}$}
\label{eq:SDPnew}
\end{equation}
Then, we readily obtain $\|G\|_{2+} \le \ogam_a \le \oogam_a$.  
Moreover, 
it is shown in \cite{Ebihara_ECC2022} that,  
for any positive filter $G_p$, the corresponding upper bound
$\ogam_a$ satisfies $\ogam_a\le \ogam\le \|G\|_2$, 
where $\ogam$ is the COP-based {\it filter-free} upper bound
obtained by replacing \rec{eq:L2+COPnew} in \rec{eq:COPnew} with \rec{eq:L2+COPnewFilterFree}.   
Similarly, for the SDP-based {\it filter-free} upper bound $\oogam$, 
we have $\oogam_a\le \oogam\le \|G\|_2$.  
Namely, by introducing any positive filter, 
we can obtain no more conservative upper bounds than
the filter-free ones.  
The matrix $Q_a$ captures the nonnegativity of the input and filtered signals 
and is referred to as the copositive multiplier.  

\subsection{Application to GAS Analysis of ReLU-Feedback Systems and 
Relationship with IQC Conditions}
\label{sec:ReLU-Feedback}

Let us revisit the GAS analysis of the nonlinear feedback system
shown in \rfig{fig:NNfb} (right) where 
$\Phi:\ \bbR^{\nw}\to \bbR^{\nw}$ is the (repeated) ReLU.  
Then, from the discussion in Subsections \ref{sub:ub} and \ref{sub:ReLU-Feedback}, 
we can conclude that the ReLU-feedback system
is GAS if $\oogam_a<1$ holds for some positive filter.  
In this way, we can use our results for the GAS analysis of
ReLU-feedback systems.  
A concise illustrative numerical example is included in 
\ref{sub:num_ub}.  

To see the relationship between our results in \ref{sub:ub} and 
IQC-based stability conditions \cite{Megretski_IEEE1997,Veenman_EJC2016,Scherer_IEEEMag2022}
in a concise fashion, let us consider the strict inequality version
of \rec{eq:L2+COPnewFilterFree} with $\gamma=1$ that is rearranged as
\begin{equation}
\scalebox{0.95}{$
\arraycolsep=1mm
\begin{array}{@{}l}
\begin{bmatrix}
PA + A^T P & PB \\
B^T P & 0
\end{bmatrix}
+
\begin{bmatrix}
 C & D \\
 0 & I_{\nw} \\
\end{bmatrix}^T
\Pi
\begin{bmatrix}
 C & D \\
 0 & I_{\nw} \\
\end{bmatrix}
\prec 0,
\end{array}$}
\label{eq:IQC}
\end{equation}
\[
\Pi= \Pi_\mathrm{COP}:=
\begin{bmatrix}
I_{\nw} & 0 \\
0 & Q-I_{\nw}
\end{bmatrix}.  
\]
On the other hand, the general IQC theory \cite{Megretski_IEEE1997,Veenman_EJC2016,Scherer_IEEEMag2022}
guarantees the GAS of the ReLU-feedback
system if there exist $P\in\bbS_{++}^n$ and $\Pi\in\bPi$ such that \rec{eq:IQC} holds
where $\bPi\subset \bbS^{2\nw}$ is \textit{the set of multipliers} defined by
\[
\scalebox{0.95}{$
\begin{array}{@{}l}
\arraycolsep=0.1mm
\bPi:=
\left\{
\Pi\in\bbS^{2\nw}:\ 
\left[
 \begin{array}{c}
  \zeta\\
  \Phi(\zeta)\\
 \end{array}
\right]^T\Pi
\left[
 \begin{array}{c}
  \zeta\\
  \Phi(\zeta)\\
 \end{array}
\right]\ge 0\ \forall \zeta \in \bbR^{\nw}
\right\}. 
\end{array}$}
\]
Here, it is obvious that $\Pi_\mathrm{COP}\in\bPi$.  
Namely, our results can be regarded as a special case of  
the IQC-based condition by employing $\Pi_\mathrm{COP}\in\bPi$.  
We can interpret that $\Pi_\mathrm{COP}\in\bPi$ has been obtained by
focusing on the $L_2$-induced norm property $\|\Phi\|\le 1$ and 
the nonnegativity property $\Phi:\bbR^{\nw}\to\bbR_+^{\nw}$ of the ReLU $\Phi$.  
Here, it is true that we can employ more sophisticated multipliers
by closely working on the ReLU-properties, 
and the latest results on this issue can be found, 
e.g., in \cite{Ebihara_ECC2024}.  

\subsection{Concrete Construction of Positive Filters}
\label{sub:filter}

As for the positive filter $G_p$ given by \rec{eq:Gp}, 
we propose to use the specific form given by
\begin{equation}
\scalebox{0.92}{$
\begin{array}{@{}l}
A_p=A_{p,\alpha,N}:=J_{\alpha,N}\otimes I_{\nw}\in\bbR^{N\nw\times N\nw},\\
B_p=B_{p,N}:=E_N\otimes I_{\nw}\in\bbR^{N\nw\times\nw}, \\
J_{\alpha,N}:=
 \left[
\renewcommand{\arraystretch}{0.5}
 \begin{array}{ccccc}
  \alpha & 1 & 0 & \cdots & 0 \\
  0 & \alpha & 1 &  \ddots & \vdots \\
  \vdots & \ddots & \ddots &  \ddots & 0 \\
  \vdots &   & \ddots &  \ddots & 1 \\
  0 & \cdots & \cdots &  0 & \alpha \\
 \end{array}
 \right]\in \bbR^{N\times N},\
E_N:=
 \left[
 \begin{array}{c}
  0 \\
  \vdots\\
  0 \\
  1
 \end{array}
 \right]\in \bbR^{N} \hspace*{-20mm} 
\end{array}$}
\label{eq:FIR}
\end{equation}
where $\alpha<0$.  
By increasing the degree $N$ of the positive filter $G_p$ given 
by \rec{eq:Gp} and \rec{eq:FIR}, 
we can construct a sequence of COPs in the form of \rec{eq:COPnew}
and SDPs in the form of \rec{eq:SDPnew}.  
In the following, we denote by 
$\ogam_{a,\alpha,N}$ and $\oogam_{a,\alpha,N}$
the optimal values of these COPs and SDPs, respectively.  
In addition, we denote by 
$A_{a,\alpha,N}$, 
$B_{a,N}$, 
$C_{a,N}$, and 
$D_{a,N}(=D)$
the coefficient matrices of the augmented system $G_a$
given by \rec{eq:Ga} corresponding to the filter of degree $N$.  
Then, regarding the effectiveness of employing higher-degree positive filters
in improving upper bounds, 
we can obtain the next result.  
\begin{theorem}
Let us consider the upper bounds of $\|G\|_{2+}$
given by $\ogam_{a,\alpha,N}$ and $\oogam_{a,\alpha,N}$ 
that are characterized respectively by 
\rec{eq:COPnew} and \rec{eq:SDPnew} with 
the positive filter $G_p$ of the form 
\rec{eq:Gp} and \rec{eq:FIR} of degree $N$.  
Then, for $N_1\le N_2$, we have
$\ogam_{a,\alpha,N_2}\le \ogam_{a,\alpha,N_1}$ and 
$\oogam_{a,\alpha,N_2}\le \oogam_{a,\alpha,N_1}$.    
\label{th:monotone}
\end{theorem}
\begin{proofof}{\rth{th:monotone}}
In the following, we prove $\ogam_{a,\alpha,N_2}\le \ogam_{a,\alpha,N_1}$.  
The proof for $\oogam_{a,\alpha,N_2}\le \oogam_{a,\alpha,N_1}$
follows similarly.  
To prove $\ogam_{a,\alpha,N_2}\le \ogam_{a,\alpha,N_1}$, 
it suffices to show that 
$\ogam_{a,\alpha,N+1}\le \ogam_{a,\alpha,N}$ holds for any $N$.  
Furthermore, this can be verified by proving that 
\rec{eq:L2+COPnew} corresponding to the filter of degree $N+1$
holds with $\gamma=\ogam_{a,\alpha,N}+\varepsilon$ for any 
$\varepsilon>0$.  

To this end, we first note from the definition of 
$\ogam_{a,\alpha,N}$ that 
for any $\varepsilon>0$ there exist 
$P_a=P_{a,\alpha,N}\in \bbS^{n+N\nw}$ and  
$Q_a=Q_{a,\alpha,N}\in \COP^{(N+1)\nw}$ such that 
\[
 \scalebox{0.90}{$
\begin{array}{@{}l}
\begin{bmatrix}
P_{a,\alpha,N}A_{a,\alpha,N} + A_{a,\alpha,N}^T P_{a,\alpha,N} & P_{a,\alpha,N}B_{a,N} \\
\ast  & - (\ogam_{a,\alpha,N}^2+2\ogam_{a,\alpha,N}\varepsilon) I_{\nw}
\end{bmatrix}\hspace*{-20mm}\\
+
\begin{bmatrix}
C_{a,N}^T \\ D_{a,N}^T
\end{bmatrix}
\begin{bmatrix}
C_{a,N}^T \\ D_{a,N}^T
\end{bmatrix}^T\!\!\!
+
\begin{bmatrix}
0_{n,(N+1)\nw} \\ I_{(N+1)\nw}
\end{bmatrix}
Q_{a,\alpha,N}
\begin{bmatrix}
0_{n,(N+1)\nw} \\ I_{(N+1)\nw}
\end{bmatrix}^T
\preceq 0.  
\end{array}$}
\]
To proceed, define
$F_N:=
\begin{bmatrix}
 I_{\nw} & 0_{\nw,(N-1)\nw} 
 \end{bmatrix}\in\bbR^{\nw\times N\nw}$.  
Then, there exist $\varepsilon_1,\varepsilon_2>0$ such that
\[
\scalebox{0.90}{$
\begin{aligned}
&\left[
  \begin{array}{ccc}
   0_{n,n}& 0 & 0 \\
   \ast & -\varepsilon_2 I_{N\nw}-\frac{\varepsilon_1}{2\alpha} F_N^TF_N & \varepsilon_2 P_pB_p\\
    \ast & \ast & -\varepsilon^2 I_{\nw}
  \end{array}
\right]\preceq 0  
\end{aligned}$}
\]  
where $P_p\in\bbS_{++}^{\np}$ is the unique solution of the 
Lyapunov equation $P_p A_p + A_p^T P_p +I_{Nnw}=0$.   
By summing up the above two inequalities and
applying the Schur complement argument by focusing on the term
$-\frac{\varepsilon_1}{2\alpha} F_N^TF_N$,  
we obtain \rec{eq:long1} given at the top of the next page.  
\begin{figure*}[t]
\begin{equation}
\scalebox{0.9}{$
\begin{array}{@{}l}
\left[
\begin{array}{ccc}
P_{a,\alpha,N}A_{a,\alpha,N} + A_{a,\alpha,N}^T P_{a,\alpha,N} +
\begin{bmatrix}
0_{n,n}  & 0 \\ 0 & -\varepsilon_2 I_{Nw} 
\end{bmatrix} 
&
P_{a,\alpha,N}B_{a,N} + \begin{bmatrix} 0_{n,n_w} \\ \varepsilon_2 P_p B_p \end{bmatrix} 
& 
\begin{bmatrix} 0_{n,n_w} \\ \varepsilon_1 F_N^T \end{bmatrix} \\ 
\ast & -(\ogam_{a,\alpha,N}+\varepsilon)^2 I_{\nw} & 0 \\
\ast & \ast & 2\varepsilon_1\alpha I_{n_w}  
\end{array}\right]
+
\begin{bmatrix}
C_{a,N}^T \\ D_{a,N}^T \\ 0_{n_w,n_z}\\ 
\end{bmatrix}
\begin{bmatrix}
C_{a,N}^T \\ D_{a,N}^T\\ 0_{n_w,n_z}
\end{bmatrix}^T\\
\hspace*{30mm}
+
\begin{bmatrix}
0_{n,(N+1)\nw} \\ I_{(N+1)\nw} \\ 0_{\nw,(N+1)\nw}
\end{bmatrix}
Q_{a,\alpha,N}
\begin{bmatrix}
0_{n,(N+1)\nw} \\ I_{(N+1)\nw} \\ 0_{\nw,(N+1)\nw}
\end{bmatrix}^T
\preceq 0.  
\end{array}$}
\label{eq:long1}
\end{equation}
\begin{equation}
\scalebox{0.90}{$
\begin{array}{@{}l}
\left[
\begin{array}{cccc}
P_{a,\alpha,N}^{11}A + A^T P_{a,\alpha,N}^{11} & 0 & 
P_{a,\alpha,N}^{12}A_{p,\alpha,N} + A^T P_{a,\alpha,N}^{12} &
P_{a,\alpha,N}^{11}B+P_{a,\alpha,N}^{12}B_p \\
 \ast &  2\varepsilon_1 \alpha I_{\nw} &
 \varepsilon_1 F_N & 0 \\\\
\ast & \ast &
 \hatP_{a,\alpha,N}^{22}A_{p,\alpha,N} + A_{p,\alpha,N}^T \hatP_{a,\alpha,N}^{22} & P_{a,\alpha,N}^{12 T}B +\hatP_{a,\alpha,N}^{22}B_p \\
\ast & \ast & 
 \ast & -(\ogam_{a,\alpha,N}+\varepsilon)^2 I_{\nw}
\end{array}\right]
+
\begin{bmatrix}
C_{a,N+1}^T \\ D_{a,N+1}^T
\end{bmatrix}
\begin{bmatrix}
C_{a,N+1}^T \\ D_{a,N+1}^T
\end{bmatrix}^T\\
\hspace*{30mm}
+
\begin{bmatrix}
0_{n,(N+2)\nw} \\ I_{(N+2)\nw}
\end{bmatrix}
\begin{bmatrix}
0_{\nw,\nw} & 0 \\ 0 & Q_{a,\alpha,N}
\end{bmatrix}
\begin{bmatrix}
0_{n,(N+2)\nw} \\ I_{(N+2)\nw}
\end{bmatrix}^T
\preceq 0.  
\end{array}$}
\label{eq:N+1}
\end{equation}
\begin{center}
=================================================================
\end{center}
\end{figure*}
For \rec{eq:long1}, by multiplying 
$
\begin{bmatrix}
I_n                 & 0_{n,N\nw}         &  0_{n,\nw}       & 0_{n,\nw} \\  
0_{\nw,n}      & 0_{\nw,N\nw}     & 0_{\nw,\nw}   & I_{\nw} \\  
0_{N\nw,n}   & I_{N\nw}             & 0_{N\nw,\nw}        &  0_{N\nw,\nw} \\  
0_{\nw,n}      & 0_{\nw,N\nw}      & I_{\nw}         & 0_{\nw,\nw} \\  
\end{bmatrix}$ from left and its transpose from right, we obtain 
\rec{eq:N+1} given at the top of the next page.  
Here, 
\[
\scalebox{0.95}{$
\begin{array}{@{}l}
\begin{bmatrix}
P_{a,\alpha,N}^{11} & P_{a,\alpha,N}^{12} \\
P_{a,\alpha,N}^{12 T} & P_{a,\alpha,N}^{22}
\end{bmatrix}:=P_{a,\alpha,N},\
P_{a,\alpha,N}^{11}\in\bbS^{n},\
P_{a,\alpha,N}^{22}\in\bbS^{N\nw},\vspace*{2mm}\\
\hatP_{a,\alpha,N}^{22}:= P_{a,\alpha,N}^{22}+\varepsilon_{2}P_p.  
\end{array}$}  
\]  
Then, \rec{eq:N+1} shows that 
\rec{eq:L2+COPnew} corresponding to the filter of degree $N+1$
holds with $\gamma=\ogam_{a,\alpha,N}+\varepsilon$ and 
\[
 \scalebox{0.95}{$
 \begin{array}{@{}l}
  P_a=P_{a,\alpha,N+1}=
   \begin{bmatrix}
    P_{a,\alpha,N}^{11} & 0 & P_{a,\alpha,N}^{12} \\
    0 & \varepsilon_1 I_{\nw} & 0 \\
    P_{a,\alpha,N}^{12 T} & 0 & \hatP_{a,\alpha,N}^{22}
   \end{bmatrix}\in \bbS^{n+(N+1)\nw},\\   
   Q_a=Q_{a,\alpha,N+1}=
   \begin{bmatrix}
    0_{\nw,\nw} & 0 \\ 0 & Q_{a,\alpha,N}
   \end{bmatrix}\in \COP^{(N+2)\nw}.  
  \end{array}$}  
\]  
This completes the proof.  
\end{proofof}
\begin{remark}
The introduction of dynamic positive filters 
in Subsections \ref{sub:ub} and \ref{sub:filter} is
motivated by the IQC-based analysis using dynamic multipliers 
\cite{Megretski_IEEE1997,Veenman_EJC2016,Scherer_IEEEMag2022}.  
Similarly to the dynamic multipliers in the IQC framework, 
the positive filters allow us to increase the 
size (freedom) of  the copositive multipliers as we can see by comparing
\rec{eq:L2+COPnew} and \rec{eq:L2+COPnewFilterFree}.  
The effectiveness of employing positive filters is
illustrated in \ref{sub:num_ub} by numerical examples.  
It is also true that the pole $\alpha<0$ in the positive filter \rec{eq:FIR}
corresponds to the pole of the basis functions
of the dynamic multipliers which should be fixed in advance.  
Since no general recipe is available on its appropriate choice,  
we usually need to carry out a line search over $\alpha<0$.  
\end{remark}
\begin{remark}
From \rth{th:monotone}, we see that 
we can construct a monotonically non-increasing sequence of 
upper bounds $\{\oogam_{a,\alpha,N}\}$ of $\|G\|_{2+}$
by increasing the degree $N$ and solving the corresponding SDP 
\rec{eq:SDPnew}.  
This is done at the expense of increased computational burden.  
Since we cannot theoretically prove the convergence of this
sequence to the true value of $\|G\|_{2+}$, 
it is quite important to obtain sharp lower bounds of $\|G\|_{2+}$
to ensure the accuracy of the computed upper bounds.  
\end{remark}
%

\section{Main Results for Uniform Infimum of $\|G\|_{p+}/\|G\|_p$}
\label{sec:uniform}

The next theorem provides our main results for Problem~2:  
analysis of the uniform infimum on the ratio of
the $L_{2+}$ induced norm to the $L_{2}$ induced norm.   
\begin{theorem}
For $p\in [1,\infty)$ and $p=\infty$, 
the uniform infimum $\nu_p^\star$ characterized by \rec{eq:nups} are given by 
\begin{equation}
 \nu_p^\star=2^{\frac{1-p}{p}}\ (p\in[1,\infty)),\quad
 \nu_\infty^\star=\frac{1}{2}.  
\label{eq:nu}
\end{equation}
Moreover, the stable LTI system $G^\star(s)=1-e^{-Ls}\ (L>0)$ 
attains
\begin{equation}
\|G^\star\|_{p+}=\nu_p^\star\|G^\star\|_{p}\ (p\in [1,\infty),\ p=\infty).   
\label{eq:Gstar}
\end{equation}
\label{th:ULBA}
\end{theorem}
%
The proof of this theorem is given in Appendix~\ref{sec:ap2}.  
We emphasize that, 
since we only rely on the linearity of
the underlying systems and the properties of the $L_p$ norms of signals in the proof, 
\rth{th:ULBA} is valid even if we extend the
set $\GLTI$ in Problem~2 to the set of linear, stable, 
{\it time-varying}, and {\it noncausal} systems including 
infinite dimensional ones.  
%

\section{Main Results for $L_{2+}$ Lower Bound Analysis Problem: Single-Input Case}
\label{sec:lbsingle}

From \rth{th:ULBA}, we see that
$
\|G\|_{2+} \ge \frac{1}{\sqrt{2}}\|G\|_2\ \forall G\in \GLTI
$.  
However, for each system $G\in\GLTI$, 
it is expected that $\|G\|_{2+}$ can be strictly larger than
$\frac{1}{\sqrt{2}}\|G\|_2$.  
In fact, the main contributions of this section are as follows:
\begin{quote}
\begin{itemize}
 \item[(i)] For any nonzero stable finite-dimensional single-input LTI system $G$, 
	    we prove $\|G\|_{2+}>\frac{1}{\sqrt{2}}\|G\|_2$.  
 \item[(ii)] For a given nonzero stable finite-dimensional single-input LTI system $G$, 
	     we provide a method to compute a lower bound that is strictly larger than 
	     $\frac{1}{\sqrt{2}}\|G\|_2$.  
\end{itemize}
\end{quote}

We derive these results 
by using basics about frequency responses of LTI systems, 
and reducing the lower bound analysis problem
into a semi-infinite programming problem \cite{Shapiro_2009}.  
The treatment of multi-input systems is deferred to the next section.  
In the following, we denote by 
$\GLTISI$ the set of stable and single-input LTI systems including infinite dimensional ones.  

\subsection{Reduction to Semi-Infinite Programming Problem}

We first recall the next very basic result. 
\begin{lemma}[\cite{Folland_2009}]
For given $a_m\in\bbR\ (m=0,\cdots,N)$, $\phi_m\in\bbR\ (m=1,\cdots,N)$,
$\omega>0$, $T=2\pi/\omega$, and $\clI=[-T/2,T/2]$, we have
\[
\scalebox{0.93}{$
\begin{aligned}
&\frac{1}{T} \int_{\clI} \left(a_0+\sum_{m=1}^N a_m \cos(m\omega t+\phi_m)\right)^2dt=a_0^2+\frac{1}{2}\sum_{m=1}^N a_m^2.  
\end{aligned}$}
\]
 \label{le:basic}
\end{lemma}

For a given $G\in\GLTISI$, let us inject the nonnegative input signal 
\begin{equation}
 w_\omega^{[N]}(t):=a_0+\cos (\omega t)+\sum_{m=2}^N a_m \cos(m\omega t)\ (t \ge 0)    
\label{eq:Norder}
\end{equation}
where we assume that $a_0,a_m\ (m=2,\cdots,N)$ are chosen such that 
$w_\omega^{[N]}(t)\ge 0\ (t\ge 0)$.  
It should be noted that, 
even though $w_\omega^{[N]}(t)\not\in L_{2+}$ and 
$w_\omega^{[N]}(t)\in L_{2\mathrm{e}+}$, 
the lower bound analysis in the following on the basis of 
$w_\omega^{[N]}(t)$ is valid since $\|G\|_{2+}$ is defined 
as the supremum over $w\in L_{2+}$ in \rec{eq:Lp+ind}.  
If we denote by $z_{\omega}^{[N]}$ the corresponding steady-state output,  
we see from the steady-state analysis of the frequency response 
and \rle{le:basic} that
\begin{equation}
\scalebox{0.90}{$
\begin{array}{@{}l}
 \|G\|_{2+}\ge 
\dfrac{\sqrt{\dfrac{1}{T}\int_{\clI} z_{\omega}^{[N]}(t)^Tz_{\omega}^{[N]}(t)dt}}{\sqrt{\dfrac{1}{T}\int_{\clI} w_\omega^{[N]}(t)^2dt}}\\
=
\dfrac{\sqrt{a_0^2|G(0)|_2^2+\dfrac{1}{2}|G(j\omega)|_2^2+\dfrac{1}{2}\displaystyle\sum_{m=2}^N a_m^2|G(jm\omega)|_2^2}}{\sqrt{a_0^2+\dfrac{1}{2}+\dfrac{1}{2}\displaystyle\sum_{m=2}^N a_m^2}}\\  
=
\dfrac{\sqrt{2a_0^2|G(0)|_2^2+|G(j\omega)|_2^2+\displaystyle\sum_{m=2}^N a_m^2|G(jm\omega)|_2^2}}{\sqrt{2a_0^2+1+\displaystyle\sum_{m=2}^N a_m^2}}.  
\end{array}$}
\label{eq:key_lb}
\end{equation}
By leaving $|G(j\omega)|_2^2$ only in the 
numerator
and taking the supremum over $\omega\in[0,\infty)$, we have
\begin{equation}
\begin{array}{@{}l}
 \|G\|_{2+}\ge  \dfrac{1}{\sqrt{2a_0^2+1+\displaystyle\sum_{m=2}^N a_m^2}} \|G\|_2.  
\end{array}
\label{eq:multigen}
\end{equation}
This result motivates us to consider the following semi-infinite programming problem:
\begin{equation}
\scalebox{0.85}{$
\begin{aligned}
& \gamma_N^\star:=\inf_{a_0,a_2,\cdots,a_N} 2a_0^2+1+\displaystyle\sum_{m=2}^N a_m^2\quad \mbox{s.t.}\vspace*{1mm}\\
& w^{[N]}(t):=a_0+\cos(t)+\sum_{m=2}^N a_m \cos(m t)\ge 0\ (\forall t \in \clI:=[-\pi,\pi]).  \hspace*{-20mm}
\end{aligned}$}
\label{eq:semiinf}
\end{equation}
For each $\gamma_N^\star$, we see 
$\|G\|_{2+}\ge \dfrac{1}{\sqrt{\gamma_N^\star}}\|G\|_2\ (\forall G\in\GLTISI)$.  
In addition, since $\gamma_N^\star$ is 
monotonically non-increasing with respect to $N\in \bbN$, 
and since $\gamma_N^\star \ge 1\ (\forall N\in\bbN)$ holds, 
the sequence $\{\gamma_N^\star\}$ converges.  
If we define $\gamma^\star:=\lim_{N\to\infty} \gamma_N^\star$, 
we readily obtain
\begin{equation}
 \|G\|_{2+}\ge   \dfrac{1}{\sqrt{\gamma^\star}}\|G\|_2\ \forall G\in \GLTISI.  
\label{eq:mainlb}
\end{equation}
%

\subsection{Effective Lower Bound Computation Methods}

From \rth{th:ULBA}, we know that 
$\|G^\star\|_{2+}=\frac{1}{\sqrt{2}}\|G^\star\|_2$
and hence $\gamma^\star\ge 2$ should hold in \rec{eq:mainlb}.  
Therefore, if we are able to construct $w(t)$ such that
\begin{equation}
\scalebox{1.0}{$
\begin{aligned}
& w(t)=a_0+\cos(t)+\sum_{m=2}^\infty a_m \cos(m t)\ge 0\ (\forall t \in \clI),\hspace*{-20mm}\\  
& 2a_0^2+1+\displaystyle\sum_{m=2}^\infty a_m^2=2,
\end{aligned}$}
\label{eq:opt_cond}
\end{equation}
then this is an optimal solution for 
the semi-infinite programming problem \rec{eq:semiinf}
in the limit case $N\to\infty$.  
With this fact in mind, let us consider the nonnegative signal 
\begin{equation}
\begin{array}{@{}l}
w^\star(t):=\max(2\cos(t),0)\ (t\ge 0)
\end{array}
\label{eq:optin}
\end{equation}
whose Fourier series expansion is given by
\begin{equation}
\scalebox{0.9}{$
\begin{aligned}
w^\star(t)=\dfrac{2}{\pi}+\cos(t)+\dfrac{4}{\pi}\sum
_{p=1}^\infty \dfrac{(-1)^{p+1}}{(2p+1)(2p-1)}\cos (2p t),   
\end{aligned}$}
\label{eq:Fourier}
\end{equation}
see, e.g., \cite{Folland_2009}.   
In \rec{eq:semiinf}, this corresponds to the case where 
\begin{equation}
\scalebox{0.9}{$
\begin{array}{@{}l}
 a_0^\star = \dfrac{2}{\pi},\ 
 a_{2p}^\star = \dfrac{4}{\pi} \dfrac{(-1)^{p+1}}{(2p+1)(2p-1)},\ 
 a_{2p+1}^\star = 0\ (p\in\bbN).  
\end{array}$}
\label{eq:opt_coef} 
\end{equation}
From Parseval's identity, we readily see that
\begin{equation}
\scalebox{0.9}{$
\begin{array}{@{}l}
2a_0^{\star 2}+1+\displaystyle\sum_{m=2}^\infty a_m^{\star 2}
=2\frac{1}{2\pi}\int_{\clI} w^\star(t)^2 dt = 2.  
\end{array}$}
\label{eq:innorm}
\end{equation}
It follows that the signal $w^\star$ given by \rec{eq:optin}
satisfies the optimality condition \rec{eq:opt_cond}
and hence is an optimal solution for 
the semi-infinite programming problem \rec{eq:semiinf}
in the limit case $N\to\infty$.  
From these results, we see that the next result holds
for any $G\in\GLTISI$:  
\begin{equation}
\scalebox{0.81}{$
 \begin{aligned}
 \|G\|_{2+}\ge \sup_{\omega>0}
\dfrac{\sqrt{2a_0^{\star 2}|G(0)|_2^2+|G(j\omega)|_2^2+\displaystyle\sum_{m=2}^\infty a_m^{\star 2}|G(jm\omega)|_2^2}}{\sqrt{2}}.  \hspace*{-10mm}
\end{aligned}$}
\label{eq:key_lb2}
\end{equation}
This expression leads us  to the next results. 
\begin{theorem}
Suppose $G\in \GLTISI\setminus\{0\}$ is finite-dimensional
where ``\hspace*{0.5mm}$0$'' stands for the zero operator.  
Then, we have $\|G\|_{2+}>\frac{1}{\sqrt{2}}\|G\|_2$.  
In particular, if $\|G\|_2=|G(0)|_2$ holds,  
then $\|G\|_{2+}=\|G\|_2$ holds.  
\label{th:nonnegativity}
\end{theorem}
\begin{proofof}{\rth{th:nonnegativity}}
We consider the following three cases:
(i) $\|G\|_2$ is attained at the 
angular frequency $\omega=0$, i.e., $\|G\|_2=|G(0)|_2$;   
(ii) $\|G\|_2$ is given as $\|G\|_2=|G(j\infty)|_2$ where $G(j\infty):=\lim_{\omega\to\infty} G(j\omega)$;   (iii) $\|G\|_2$ is attained at $\omega=\omega^\star\in (0,\infty)$, i.e., 
$\|G\|_2=|G(j\omega^\star)|_2$.  

(i) Suppose $\|G\|_2=|G(0)|_2$.  
Then, by letting $\omega\to 0$
in \rec{eq:key_lb2}, we see from \rec{eq:innorm} that 
$\|G\|_{2+}\ge |G(0)|_2=\|G\|_2$.  Namely, $\|G\|_{2+}=\|G\|_2$ holds.  

(ii) Suppose $\|G\|_2=|G(j\infty)|_2$.  
Then, we see from \rec{eq:key_lb2} and \rec{eq:innorm} that
\begin{equation}
\scalebox{0.87}{$
\begin{aligned}
\|G\|_{2+}
\ge& \frac{\sqrt{2a_0^{\star 2}|G(0)|_2^2+(2-2a_0^{\star 2})|G(j\infty)|_2^2}}{\sqrt{2}}\\
\ge& \frac{\sqrt{(2-2a_0^{\star 2})}}{\sqrt{2}}|G(j\infty)|_2
\approx \frac{1.0906}{\sqrt{2}}|G(j\infty)|_2
>\frac{1}{\sqrt{2}}\|G\|_2.\hspace*{-20mm}
\end{aligned}$}
\label{eq:evalinf}
\end{equation}

(iii) Suppose $\|G\|_2=|G(j\omega^\star)|_2\ (\omega^\star\in(0,\infty))$.  
Then, for $\|G\|_{2+}=\frac{1}{\sqrt{2}}\|G\|_2$ to hold, we see 
from \rec{eq:opt_coef} and \rec{eq:key_lb2} that the system $G$ should satisfy 
the infinitely many interpolation constraints:
\begin{equation}
 G(0)=0,\ G(j2p\omega^\star)=0\ (p\in\bbN).  
\label{eq:interpolation}
\end{equation}
This is impossible for the finite-dimensional system 
$G\in \GLTISI\setminus\{0\}$ and hence $\|G\|_{2+}>\frac{1}{\sqrt{2}}\|G\|_2$.  
\end{proofof}
\begin{remark}
In \rth{th:ULBA}, we have shown that 
the infinite-dimensional system $G^\star(s)=1-e^{-Ls}\ (L>0)$ 
satisfies $\|G^\star\|_{2+}=\frac{1}{\sqrt{2}}\|G^\star\|_{2}$.  
Therefore, from the proof of \rth{th:nonnegativity}, 
the system $G^\star$ should satisfy the infinitely many interpolation constraints 
\rec{eq:interpolation}.  Indeed, we see that $\omega^\star=\frac{1}{L}\pi$
for $G^\star$, and $G^\star$ 
does satisfy the interpolation constraints \rec{eq:interpolation} since
$G^\star(0)=0$, and $G^\star(j2p\omega^\star)=0\ (p\in\bbN)$.  
\end{remark}

For a given system $G\in\GLTISI$, 
we finally make active use of \rec{eq:key_lb2} 
for the lower bound computation of $\|G\|_{2+}$.  
By truncation of the infinite series, let us define 
\begin{equation}
\scalebox{0.83}{$
\begin{aligned}
& \upsilon_N(G):=\sup_{\omega>0}
\dfrac{\sqrt{2a_0^{\star 2}|G(0)|_2^2+|G(j\omega)|_2^2+\displaystyle\sum_{m=2}^N a_m^{\star 2}|G(jm\omega)|_2^2}}{\sqrt{2}}.  \hspace*{-10mm}
\end{aligned}$}
\label{eq:upsilon_N}
\end{equation}
Then, it is straightforward from \rth{th:nonnegativity} that the next results hold.  
\begin{theorem}
Suppose $G\in \GLTISI\setminus\{0\}$ is finite-dimensional and define
$\upsilon_N(G)\ (N\in\bbN)$ by \rec{eq:upsilon_N}.  
Then, we have
\[
 \|G\|_{2+}\ge \upsilon_N(G)\ge \frac{1}{\sqrt{2}}\|G\|_2\ (\forall N\in \bbN).  
\]
In particular, $\upsilon_N(G)$ is monotonically non-decreasing with respect to $N\in\bbN$, 
and for sufficiently large $N$ we have
$\upsilon_N(G)> \frac{1}{\sqrt{2}}\|G\|_2$.  
\label{th:lbcomp} 
\end{theorem}
We finally note that $\upsilon_N(G)$ can readily be computed 
by using \rec{eq:opt_coef} in the way that
\begin{equation}
\scalebox{0.90}{$
\begin{array}{@{}l}
 \upsilon_N(G)=\|\hatG_N\|_2,\ 
 \hatG_N(s):=
\frac{1}{\sqrt{2}}
\left[
\begin{array}{c}
\sqrt{2}a_0^\star G(0)\\
G(s)\\
a_2^\star G(2s)\\
\vdots\\
a_N^\star G(Ns)\\
\end{array}
\right].  
\end{array}$}
\label{eq:upsilon_comp}
\end{equation}
%

\section{Main Results for $L_{2+}$ Lower Bound Analysis Problem: Multi-Input Case}
\label{sec:lbmulti}

To see the difference between single-input and multi-input  cases, 
let us consider the ``static'' multi-input system $G:=[\ 1\ -1\ ]$.  
Then, we can readily see that
$\|G\|_2=\sqrt{2}$ and $\|G\|_{2+}=1$.  
Namely, $\|G\|_{2+}=\frac{1}{\sqrt{2}} \|G\|_{2}$ holds.  
This example clearly shows that whole assertions in 
Theorems \ref{th:nonnegativity} and \ref{th:lbcomp} do not 
hold for multi-input systems in general.  
With this fact in mind, the goal of this section is as follows:  
\begin{enumerate}
 \item[(i)] For any $G\in \GLTI$, 
	    we provide a method to compute a lower bound that is {\it larger than or equal to}
	    $\frac{1}{\sqrt{2}}\|G\|_2$.  
 \item[(ii)] By numerical examples on a multi-input system $G\in \GLTI$, 
	     we show that the method can yield lower bounds
	     that are strictly larger than $\frac{1}{\sqrt{2}}\|G\|_2$.  
\end{enumerate}
On the size of the input and output signals in \rec{eq:G},  
we remind $w(t)\in\bbR^{\nw}$ and $z(t)\in\bbR^{\nz}$.

\subsection{The Case $\|G\|_2=\|G(j\omega^\star)\|_2\ (\omega^\star\in(0,\infty))$}

We first consider the case where $\|G\|_2=\|G(j\omega^\star)\|_2\ (\omega^\star\in(0,\infty))$ 
holds, 
i.e., $\|G\|_2$ is attained at the angular frequency $\omega=\omega^\star\ (\omega^\star\in(0,\infty))$.  
Suppose $v\in\bbC^{\nw}$ is 
the unit right singular vector corresponding to 
the maximum singular value of $G(j\omega^\star)$.  
In the following, we represent the $i$-th element of $v$ by
$v_i=|v_i|e^{j\theta_i}\ (i=1,\cdots,\nw)$.  

Under these notations, let us consider the 
nonnegative input signal $w^\star(t)\in\bbR^{\nw}$ defined by
\begin{equation}
w_i^\star(t):=|v_i|\times \max(2\cos (\omega^\star t+\theta_i),0) \ (i=1,\cdots,\nw).  
\label{eq:wstar_multi} 
\end{equation}
If we define $v^{[m]}\in\bbC^{\nw}\ (m=0,1,2\cdots)$ by 
\begin{equation}
v_i^{[m]} =|v_i|e^{jm\theta_i}\ (i=1,\cdots,\nw),   
\label{eq:vs}
\end{equation}
we see from \rec{eq:optin}, \rec{eq:Fourier}, \rec{eq:opt_coef} that
\begin{equation}
\begin{aligned}
& w^\star(t) =\Re(\hatw^\star(t)),\\ 
& \hatw^\star(t)=a_0^\star v^{[0]} + v^{[1]}e^{j\omega^\star t}+
\sum_{m=2}^\infty a_m^\star v^{[m]}e^{jm\omega^\star t}.  
\end{aligned}
\label{eq:hatw_star}
\end{equation}
Then, if we define single-input transfer functions by
\begin{equation}
G^{[m]}(s) :=G(s)v^{[m]}\ (m=0,1,2\cdots),   
\label{eq:Gs}
\end{equation}
we see from the basic property of frequency response of LTI systems that
the steady state output $z^\star(t)\in\bbR^{\nz}$, 
corresponding to the input $w^\star(t)\in\bbR^{\nw}$, is given by
\[
 \scalebox{0.85}{$
\begin{aligned}
& z^\star(t) =\Re(\hatz^\star(t)),\\ 
& \hatz^\star(t)=a_0^\star G^{[0]}(0) + G^{[1]}(j\omega^\star)e^{j\omega^\star t}+
\sum_{m=2}^\infty a_m^\star G^{[m]}(jm\omega^\star)e^{jm\omega^\star t}.  
\end{aligned}$}
\]
We also note that 
$\|G\|_2=\|G(j\omega^\star)\|_2=|G^{[1]}(j\omega^\star)|_2$.  
Then, for $T=2\pi/\omega^\star$ and $\clI=[-T/2,T/2]$, we have
\begin{equation}
\scalebox{0.80}{$
\begin{aligned}
& \|G\|_{2+}\ge 
\dfrac{\sqrt{\dfrac{1}{T}\int_{\clI} z^\star(t)^T z^\star(t)dt}}{\sqrt{\dfrac{1}{T}\int_{\clI} 
w^\star(t)^Tw^\star(t) dt}}\\
& =
\dfrac{\sqrt{2a_0^{\star 2}|G^{[0]}(0)|_2^2+|G^{[1]}(j\omega^\star)|_2^2+\displaystyle\sum_{m=2}^\infty a_m^{\star 2}|G^{[m]}(jm\omega^\star)|_2^2}}
{\sqrt{2}}  \\
& \ge \dfrac{1}{\sqrt{2}}|G^{[1]}(j\omega^\star)|_2  = \dfrac{1}{\sqrt{2}}\|G\|_2.  
\end{aligned}$}
\label{eq:multi}
\end{equation}
Motivated by this fact, let us define
\begin{equation}
\scalebox{0.72}{$
\begin{aligned}
& \upsilon_N(G):=\sup_{\omega>0}
\dfrac{\sqrt{2a_0^{\star 2}|G^{[0]}(0)|_2^2+|G^{[1]}(j\omega)|_2^2+\displaystyle\sum_{m=2}^N a_m^{\star 2}|G^{[m]}(jm\omega)|_2^2}}{\sqrt{2}}.  \hspace*{-10mm}
\end{aligned}$}
\label{eq:upsilon_N_multi}
\end{equation}
Then, it is straightforward from \rec{eq:multi} that the next results hold.  
\begin{theorem}
Suppose $G\in \GLTI$ and define
$\upsilon_N(G)\ (N\in\bbN)$ by \rec{eq:upsilon_N_multi}.  
Then, we have
\[
 \|G\|_{2+}\ge \upsilon_N(G)\ge \frac{1}{\sqrt{2}}\|G\|_2\ (\forall N\in \bbN).  
\]
In particular, $\upsilon_N(G)$ is monotonically non-decreasing with respect to $N\in\bbN$.   
\label{th:lbcomp_multi} 
\end{theorem}
%
%
\begin{remark}
In the single-input case (i.e., if $\nw=1$), it is clear that
\rth{th:lbcomp_multi} reduces to \rth{th:lbcomp}.  
This is because if $\nw=1$ then we have
$v^{[m]}=1\ (m=0,1,2\cdots)$ and hence
$G^{[m]}(s)=G(s)\ (m=0,1,2\cdots)$ hold in \rec{eq:Gs}.  
\end{remark}
%

\subsection{The Case $\|G\|_2=\|G(0)\|_2$}

We next consider the case where $\|G\|_2=\|G(0)\|_2$ holds, 
i.e., $\|G\|_2$ is attained at the angular frequency $\omega=0$.  
We cannot simply apply the method in 
the preceding subsection to $\omega^\star=0$, 
since our underlying result, \rle{le:basic}, never holds for $\omega=0$.  
Similarly to the single-input case, 
such a discussion to take the limit $\omega^\star\to 0$ might be plausible, 
but we encounter a difficulty since, 
as shown in \rec{eq:vs}, 
the method in the preceding subsection 
requires an explicit representation of the unit right singular vector corresponding to 
the maximum singular value of $G(j\omega^\star)$ (this is not the case for single-input case).  
On the basis of these observations, here 
we first blindly follow the discussion in the preceding subsection
applied to $\omega^\star=0$, and validate the obtained results directly later on.  

Suppose $v\in\bbR^{\nw}$ is 
the unit right singular vector corresponding to 
the maximum singular value of the matrix $G(0)\in\bbR^{\nz\times \nw}$.  
We define $v_+=\max(v,0)$ and $v_-=-\max(-v,0)$.  
Then, we see that $w^\star$ defined by \rec{eq:wstar_multi}, 
corresponding to $\omega^\star=0$, becomes 
$w^\star(t) = 2v_+$
where $\theta_i\ (i=1,\cdots,\nw)$ in \rec{eq:wstar_multi} are $0$ or $\pi$.  
We also note that $w^\star$ given above can be rewritten as 
$w^\star(t) = v+\vabs$
where $v_{\mathrm{abs},i}=|v_i|\ (i=1,\cdots,\nw)$.   
On the other hand, we see that $v^{[m]}\ (m=0,1,2,\cdots)$
defined by \rec{eq:vs} become
\begin{equation}
 v^{[m]}=\left\{
\begin{array}{ll}
 v & m\ \mbox{is odd},\\
 \vabs & m\ \mbox{is even}.  \\
\end{array}
 \right.  
\label{eq:vs2}
\end{equation}
Then,  the inequality \rec{eq:multi} reduces to
\[
\scalebox{0.90}{$
\begin{aligned}
\|G\|_{2+} & \ge 
\dfrac{\sqrt{2a_0^{\star 2}|G(0)v^{[0]}|_2^2+|G(0)v^{[1]}|_2^2+\displaystyle\sum_{m=2}^\infty a_m^{\star 2}|G(0)v^{[m]}|_2^2}}
{\sqrt{2}}  \\
& =  \dfrac{1}{\sqrt{2}}\sqrt{\|G(0)\|_2^2+|G(0)\vabs|_2^2}  \\
\end{aligned}$}
\]
where we used \rec{eq:opt_coef}, \rec{eq:innorm}, and \rec{eq:vs2}.  
Since $\|G\|_{2+}\ge \|G(0)\|_{2+}$ obviously holds, 
it follows that if we are able to directly validate
\[
\|G(0)\|_{2+}\ge 
\dfrac{1}{\sqrt{2}}\sqrt{\|G(0)\|_2^2+|G(0)\vabs|_2^2} 
\]
for $G(0)\in\bbR^{\nz\times\nw}$, 
then we arrive at the conclusion that 
all of the results in the preceding subsection are valid even for $\omega^\star=0$.  
This validation is indeed possible as we see in the next theorem.  
\begin{theorem}
For a given $M\in\bbR^{n\times m}$, suppose $v\in\bbR^m$ is 
the unit right singular vector corresponding to the maximal singular value 
of $M$.  
Let us define $\vabs\in\bbR^{m}_+$ 
by $v_{\mathrm{abs},i}=|v_i|\ (i=1,\cdots,m)$.  
Then, we have
\begin{equation}
\|M\|_{2+}\ge \frac{1}{\sqrt{2}}\sqrt{\|M\|_2^2+|M \vabs|_2^2}.  
\label{eq:lb_Mat} 
\end{equation}
\label{th:lb_Mat} 
\end{theorem}
\begin{proofof}{\rth{th:lb_Mat}}
Without loss of generality, 
we choose (the sign of) $v$ such that
$|v_+|_2\ge |v_-|_2$ where $v_+=\max(v,0)$ and $v_-=-\max(-v,0)$.  
To prove \rec{eq:lb_Mat}, it suffices to show that
\[
 \dfrac{|M v_+|_2}{|v_+|_2}=\dfrac{|M(v+\vabs)|_2}{|v+\vabs|_2}\ge \frac{1}{\sqrt{2}}\sqrt{\|M\|_2^2+|M \vabs|_2^2}
\]
or equivalently, 
\[
|M(v+\vabs)|_2\ge \frac{1}{\sqrt{2}}|v+\vabs|_2 \sqrt{\|M\|_2^2+|M \vabs|_2^2}.  
\]
Since $v^T\vabs=|v_+|_2^2-|v_-|_2^2(\ge 0)$ holds, we obtain
\[
|M(v+\vabs)|_2^2=\|M\|_2^2(1+2(|v_+|^2-|v_-|^2))+|M\vabs|_2^2.  
\]
On the other hand, we have
\[
\begin{aligned}
& \left(\frac{1}{\sqrt{2}}|v+\vabs|_2 \sqrt{\|M\|_2^2+|M \vabs|_2^2}\right)^2\\
& =\frac{1}{2}|v+\vabs|_2^2(\|M\|_2^2+|M\vabs|_2^2)\\
& =(1+|v_+|^2-|v_-|^2)(\|M\|_2^2+|M\vabs|_2^2).  
\end{aligned}
\]
Therefore we see that
\[
\begin{aligned}
& |M(v+\vabs)|_2^2- \left(\frac{1}{\sqrt{2}}|v+\vabs|_2 \sqrt{\|M\|_2^2+|M \vabs|_2^2}\right)^2\\
& =(|v_+|^2-|v_-|^2)\|M\|_2^2-(|v_+|^2-|v_-|^2)|M\vabs|_2^2\\
& =(|v_+|^2-|v_-|^2)(\|M\|_2^2-|M\vabs|_2^2)\ge 0.  
\end{aligned}
\]
This completes the proof.  
\end{proofof}
%
%
\begin{remark}
From the proof of \rth{th:lb_Mat}, we see that
\[
 \|M\|_{2+}\ge \frac{|Mv_+|_2}{|v_+|_2}\ge \frac{1}{\sqrt{2}}\|M\|_2.  
\]
Therefore, we can obtain a lower bound of $\|M\|_{2+}$ that is larger (no smaller) than
$\frac{1}{\sqrt{2}}\|M\|_2$ by computing $|Mv_+|_2/|v_+|_2$.  
\end{remark}
%

\subsection{The Case $\|G\|_2=\|G(\infty)\|_2$}

We finally consider the case $\|G\|_2=\|G(j\infty)\|_2$
where $G(j\infty):=\lim_{\omega\to\infty} G(j\omega)$.  
We assume that $G(j\infty)\in\bbR^{\nz\times \nw}$ is well-defined; 
this is typically the case where $G$ is a finite-dimensional LTI system of the form 
\rec{eq:GLTI} and in this case we have $G(j\infty)=D\in\bbR^{\nz\times \nw}$.  

Suppose $v\in\bbR^{\nw}$ is the unit right singular vector corresponding to 
the maximum singular value of $G(j\infty)\in\bbR^{\nz\times \nw}$.  
We define $v^{[m]}=v\ (m=0,1,2,\cdots)$ and define 
$G^{[m]}(s)\ (m=0,1,2,\cdots)$ by \rec{eq:Gs}.  
Then, it is not hard to see that all the results in \rth{th:lbcomp_multi} hold.  
In particular, similarly to \rec{eq:evalinf}, we have
\[
\scalebox{0.80}{$
\begin{aligned}
 \|G\|_{2+}\ge &
\frac{\sqrt{2a_0^{\star 2}|G(0)v|_2^2+(2-2a_0^{\star 2})|G(j\infty)v|_2^2}}{\sqrt{2}}\\
\ge & \frac{\sqrt{(2-2a_0^{\star 2})}}{\sqrt{2}}|G(j\infty)v|_2
\approx \frac{1.0906}{\sqrt{2}}|G(j\infty)v|_2> \frac{1}{\sqrt{2}}\|G\||_2.  
\end{aligned}$}
\]
%

\subsection{Summary of the Lower Bound Computation}
\label{sub:lb_summary}

For a given $G\in\GLTI$, 
we now summarize the procedure for the lower bound computation 
of $\|G\|_{2+}$.  
\begin{itemize}
 \item[1.] If $\|G\|_2=\|G(j\omega^\star)\|\ (\omega^\star\in[0,\infty))$, 
	   denote by $v\in\bbC^{\nw}$ ($v\in\bbR^{\nw}$ if $\omega^\star=0$)
	   the unit right singular vector corresponding to the
	   maximal singular value of $G(j\omega^\star)$.  
	   Define $v^{[m]}\ (m=0,1,2\cdots)$ by \rec{eq:vs}.  
	   If $\|G\|_2=\|G(j\infty)\|$, 
	   denote by $v\in\bbR^{\nw}$ the unit right singular vector corresponding to 
	   the maximal singular value of $G(j\infty)\in\bbR^{\nz\times \nw}$.  
	   Define $v^{[m]}=v\ (m=0,1,2\cdots)$.  
 \item[2.] Define $G^{[m]}(s)\ (m=0,1,2\cdots)$ by \rec{eq:Gs} and compute
	   $\upsilon_N(G)$ defined by \rec{eq:upsilon_N_multi}.  
\end{itemize}
Then, as stated in \rth{th:lbcomp_multi} , we have 
\[
 \|G\|_{2+}\ge \upsilon_N(G)\ge \frac{1}{\sqrt{2}}\|G\|_2\ (\forall N\in \bbN).  
\]
In particular, $\upsilon_N(G)$ is monotonically non-decreasing with respect to $N\in\bbN$.   

\begin{remark}
Since $v\in\bbC^{\nw}$ in general, the transfer functions 
$G^{[m]}(s)\ (m=0,1,2\cdots)$ in \rec{eq:Gs} have complex coefficients.  
Still, the computation of $\upsilon_N(G)$ in \rec{eq:upsilon_N_multi} can be done
similarly to \rec{eq:upsilon_comp}.  
\end{remark}
%

\subsection{Numerical Examples}
\label{sub:num}

\subsubsection{$L_{2+}$ Upper and Lower Bounds Computation }
\label{sub:num_ub}

\begin{figure}[t]
\begin{center}
\includegraphics[width=8.0cm]{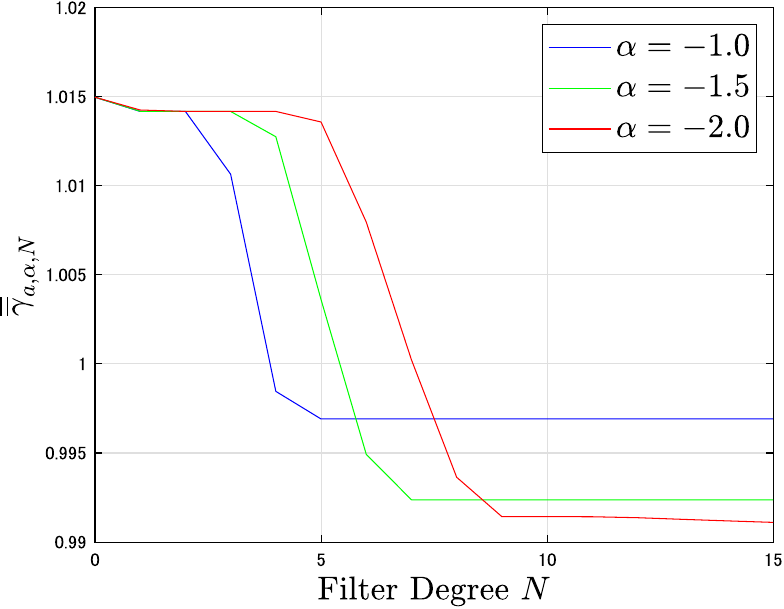}
\caption{Computed Upper Bounds $\oogam_{a,\alpha,N}$.  } 
\label{fig:ub_MI}
\end{center}
\begin{center}
\includegraphics[width=7.85cm]{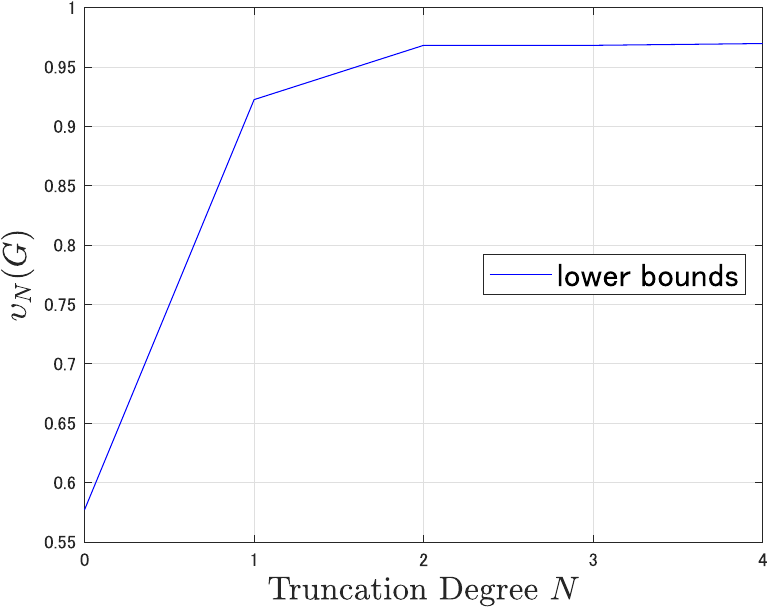}
\caption{Computed Lower Bounds $\upsilon_N(G)$.  } 
\label{fig:lb_MI}
\end{center}
\end{figure}

Let us consider the case where the system $G$ in \rec{eq:G} 
is given by \rec{eq:GLTI} with
\[
\scalebox{0.73}{$
\arraycolsep=0.5mm
\begin{array}{@{}l}
 A= 
\begin{bmatrix*}[r]
 -0.88 & 0.32 & 0.72 & -0.39 & -0.45 & -0.23 \\
  0.00 & -1.75 & -0.75 & -0.91 & 1.00 & 0.64 \\
  0.06 & -0.56 & -1.40 & 0.67 & 0.67 & 0.19 \\
  0.82 & 0.21 & -0.13 & -1.95 & 0.58 & 0.56 \\
  0.16 & -0.63 & 0.23 & 0.32 & -1.38 & 0.98 \\
  0.55 & -0.60 & 0.47 & 0.79 & 0.09 & -1.11 \\
 \end{bmatrix*},\ 
 B=
 \begin{bmatrix*}[r]
  0.32 & -0.44 & 0.14 \\
 -0.01 & 0.30 & -0.23 \\
  0.34 & -0.21 & -0.02 \\
 -0.16 & 0.17 & 0.40 \\
 -0.27 & 0.27 & -0.23 \\
  0.00 & 0.17 & -0.27 \\
 \end{bmatrix*}, \vspace*{1mm}\\
 C= 
  \begin{bmatrix*}[r]
 -0.23 & 0.10 & -0.19 & 0.43 & -0.45 & 0.33 \\
  0.27 & -0.31 & -0.27 & -0.46 & -0.07 & 0.12 \\
 -0.17 & -0.41 & 0.40 & 0.09 & 0.02 & -0.24 \\
  \end{bmatrix*},\ 
  D=
  \begin{bmatrix*}[r]
  0.61 & 0.59 & 0.18 \\
  0.39 & -0.20 & -0.62 \\
  0.10 & -0.24 & -0.19 \\
  \end{bmatrix*}. 
\end{array}$}
\]
These matrices were randomly generated. 
In this case, it turned out that $\|G\|_2=\|G(j\omega^\star)\|_2=1.0178$ where
$\omega^\star=0.6529$.  
We first computed upper bounds of $\|G\|_{2+}$
by following the method in \rsec{sec:ub}.  
The filter-free upper bound turned out to be $1.0150$.  
Then, for the poles of the positive filters, we tested the three cases:
$\alpha\in\{-1.0,-1.5,-2.0\}$.  The results are shown in \rfig{fig:ub_MI}.  
The best (least) upper bound with positive filters was $0.9911\ (\alpha=-2.0,\ N=15)$.  
With these facts in mind, we next computed $\upsilon_N(G)$, 
the lower bounds of $\|G\|_{2+}$, 
by following the procedure summarized in the preceding subsection.  
The results are shown in \rfig{fig:lb_MI}.  
The best (largest) lower bound is $0.9698$.  
This is indeed larger than the lower bound
$\frac{1}{\sqrt{2}}\|G\|_2\approx 0.7197$
obtained from the uniform infimum shown in \rth{th:ULBA}.  
From these upper/lower bounds together with
$(0.9911-0.9698)/0.9911\approx 0.0215$, 
we can conclude that the relative error between 
the upper bound $0.9911$ and the true value of $\|G\|_{2+}$
is less than $2.2\%$.  

In \rfig{fig:lb_MI}, the improvement of 
$\upsilon_N(G)$ over $N\ge 2$ is not significant.  
This happens in general due to the following reasons.  
For $N\ge 2$, we have
from \rec{eq:upsilon_N_multi} and \rec{eq:innorm} that
\[
\scalebox{0.85}{$
\begin{aligned}
\upsilon_N(G)^2
&=\sup_{\omega>0}
\dfrac{2a_0^{\star 2}|G^{[0]}(0)|_2^2+|G^{[1]}(j\omega)|_2^2+\displaystyle\sum_{m=2}^N a_m^{\star 2}|G^{[m]}(jm\omega)|_2^2}{2}\\
&\le \sup_{\omega>0}
\dfrac{2a_0^{\star 2}|G^{[0]}(0)|_2^2+|G^{[1]}(j\omega)|_2^2+(1-2a_0^{\star 2})\|G\|_2^2}{2}\\
& =\upsilon_1(G)^2+\dfrac{(1-2a_0^{\star 2})\|G\|_2^2}{2}.  
\end{aligned}$}
\]
Therefore, for $N\ge 2$, we have
\[
\dfrac{\upsilon_N(G)^2-\upsilon_1(G)^2}{\|G\|_2^2}\le \dfrac{(1-2a_0^{\star 2})}{2}\approx 0.0947.  
\]
This clearly shows that the improvement of $\upsilon_N(G)$ over $N\ge 2$ is not significant in general.  

Finally, suppose $G$ in this numerical example is the LTI
system of the ReLU-feedback system shown in \rfig{fig:NNfb} (right), and
let us illustrate how we can use lower/upper bounds of $\|G\|_{2+}$
for the GAS analysis.  
First, since $\|G\|=1.0178>1$, we cannot conclude the GAS
from the standard $L_{2}$-induced-norm-based small-gain theorem.  
On the other hand, since an easy-to-compute lower bound of $\|G\|_{2+}$
from \rth{th:ULBA}
turns out to be  
$\frac{1}{\sqrt{2}}\|G\|_2\approx 0.7197$, 
it remains possible to have $\|G\|_{2+}<1$
(if $\frac{1}{\sqrt{2}}\|G\|_2>1$ then we have
no need to resort to the upper bound computation
of $\|G\|_{2+}$ since $\|G\|_{2+}<1$ never happens).  
Finally, since it turns out that 
$\|G\|_{2+}<0.9911<1$ by the upper bound computation, 
we can conclude the GAS of the corresponding ReLU-feedback system 
using the $L_{2+}$-induced-norm-based small-gain theorem.  

%

\subsubsection{Difference of Positive Systems}

Let us consider stable and internally (and hence externally) positive systems
$G_1$, $G_2$, $G_3$ of the form \rec{eq:GLTI}
whose coefficient matrices are
\[
\scalebox{0.75}{$
\arraycolsep=0.5mm
\begin{array}{@{}l}
 A_1= 
\begin{bmatrix*}[r]
   -2.67 &   0.44 &   0.23 &   0.11 &   0.98 &   0.95 \\ 
    0.94 &  -2.93 &   0.64 &   0.53 &   0.95 &   0.11 \\ 
    0.48 &   0.98 &  -2.81 &   0.21 &   0.47 &   0.12 \\ 
    0.40 &   0.44 &   0.12 &  -3.04 &   0.99 &   0.00 \\ 
    0.19 &   0.80 &   0.41 &   0.99 &  -2.22 &   0.30 \\ 
    0.34 &   0.82 &   0.07 &   0.74 &   0.31 &  -2.70 \\ 
 \end{bmatrix*},\ 
 B_1=
 \begin{bmatrix*}[r]
    0.12 &   0.11 \\ 
    0.71 &   0.72 \\ 
    0.96 &   0.88 \\ 
    0.32 &   0.00 \\ 
    0.79 &   0.10 \\ 
    0.67 &   0.37 \\ 
 \end{bmatrix*}, \\
 C_1= 
  \begin{bmatrix*}[r]
    \ 0.91 &   \ 0.17 &   \ 0.10 &   \ 0.41 &   \ 0.20 &   \ 0.88 \\ 
    1.00 &   0.17 &   0.61 &   0.26 &   0.19 &   0.63 \\ 
  \end{bmatrix*},\ 
  D_1=
  \begin{bmatrix*}[r]
    0.62 &   0.18 \\ 
    0.66 &   0.63 \\ 
  \end{bmatrix*}, \\
 A_2= 
\begin{bmatrix*}[r]
   -1.43 &   0.70 &   0.04 &   0.70 \\
    0.49 &  -1.28 &   0.49 &   0.32 \\
    0.99 &   0.50 &  -1.20 &   0.78 \\
    0.39 &   0.06 &   0.06 &  -1.31 \\
 \end{bmatrix*},\ 
 B_2=
 \begin{bmatrix*}[r]
    0.02 &   0.25 \\
    0.78 &   0.37 \\
    0.41 &   0.84 \\
    0.48 &   0.88 \\
 \end{bmatrix*}, \\
 C_2= 
  \begin{bmatrix*}[r]
    0.24 &   0.27 &   0.10 &   0.96 \\
    0.70 &   0.44 &   0.30 &   0.42 \\
  \end{bmatrix*},\ 
  D_2=
  \begin{bmatrix*}[r]
    0.92 &   0.70 \\
    0.63 &   0.55 \\
  \end{bmatrix*}, \\
 A_3= 
\begin{bmatrix*}[r]
   -1.19 &   0.33 &   0.06 &   0.17 \\
    0.92 &  -1.57 &   0.90 &   0.42 \\
    0.26 &   0.74 &  -1.32 &   0.00 \\
    0.58 &   0.67 &   0.44 &  -1.50 \\
 \end{bmatrix*},\ 
 B_3=
 \begin{bmatrix*}[r]
    0.56 &   0.87 \\
    0.07 &   0.75 \\
    0.39 &   0.75 \\
    0.50 &   0.85 \\
 \end{bmatrix*}, \\
 C_3= 
  \begin{bmatrix*}[r]
    0.24 &   0.96 &   0.10 &   0.12 \\
    0.34 &   0.14 &   0.91 &   0.74 \\
  \end{bmatrix*},\ 
  D_3=
  \begin{bmatrix*}[r]
    0.72 &   0.98 \\
    0.12 &   0.21 \\
  \end{bmatrix*}.   
\end{array}$}
\]
These matrices were randomly generated under positivity constraints. 
Here, we regard $G_2$ and $G_3$ of dimension four 
as reduced order models of $G_1$ of dimension six.  
It turned out that $\|G_1-G_2\|_2=12.43$ and $\|G_1-G_3\|_2=15.69$
and hence we could say that $G_2$ is better than $G_3$ as a reduced order model
in the standard $L_2$ induced norm sense.  
However, since $G_1$, $G_2$, and $G_3$ are all positive, 
it is more appropriate to compare the closeness in terms of 
$L_{2+}$ induced norm as stated in Subsection \ref{sub:diff_pos}.  
By the upper bound computation method in \rsec{sec:ub}  and
the lower bound computation method in \rsub{sub:lb_summary}, 
we obtained
$12.31\le \|G_1-G_2\|_{2+}\le 12.37$ and 
$11.23\le \|G_1-G_3\|_{2+}\le 11.89$.  
Therefore we are led to the definite conclusion that 
$G_3$ is better than $G_2$ as a reduced order model since
$\|G_1-G_3\|_{2+}<\|G_1-G_2\|_{2+}$.  
This definite conclusion cannot be obtained unless the proposed
upper and lower bounds computation methods.

\section{Conclusion}

In this paper, we introduced $L_{p+}\ (p\in[1,\infty),\ p=\infty)$ induced norms for 
continuous time LTI systems.  
We first developed an $L_{2+}$ upper bound computation method 
(Theorems~\ref{th:main1} and \ref{th:monotone}), 
we then derived 
in an explicit fashion 
the uniform infimum on the ratio of
the $L_{2+}$ induced norm to the $L_{2}$ induced norm 
over all linear systems including infinite-dimensional ones 
 (Theorem~\ref{th:ULBA}), 
and finally derived an effective method to compute
lower bounds of the $L_{2+}$ induced norm 
that are better (no smaller) than the value determined 
by the uniform infimum (Theorems~\ref{th:lbcomp} and \ref{th:lbcomp_multi}).  

Even though we illustrated  the effectiveness of 
the proposed upper/lower bounds computation methods by numerical examples, 
it is true that there is no theoretical guarantees for the accuracy of these bounds.  
It is our important future topic 
to enable the quantitative evaluation of the accuracy of these bounds.   

On the other hand, as stated in Subsection~\ref{sec:ReLU-Feedback}, 
from the IQC framework 
\cite{Megretski_IEEE1997,Veenman_EJC2016,Scherer_IEEEMag2022},  
the COP for the upper bound computation of the $L_{2+}$ induced norm 
can be regarded as a result
that employs copositive multipliers capturing the nonnegativity of the assumed nonlinearities.  
It is also our future topic to clarify the effectiveness of the copositive multipliers
in static and dynamical nonlinear system analysis from a broad perspective.  


\appendix

\section{Proof of \rpr{pr:expos_result}}
\label{sec:ap0}
%
\begin{proofof}{\rpr{pr:expos_result}}
Suppose $G$ is externally positive.  
For $w\in L_p$ with $\|w\|_p=1$, 
let us define $w_+,w_-\in L_{p+}$ such that $w=w_+-w_-$ by 
$w_+(t)=\max(w(t),0_n)$, 
$w_-(t)=\max(-w(t),0_n)\ (t\in[0,\infty))$.  
Then, for $p\in[1,\infty)$ and $p=\infty$, we have
\[
 \begin{array}{@{}lcl}
 \|Gw\|_p&=& \|G(w_+-w_-)\|_p\\
 &=& \|Gw_+-Gw_-\|_p\\
 &\le & \|Gw_+ + Gw_-\|_p\\
 &= & \|G(w_+ + w_-)\|_p\\
 &\le & \|G\|_{p+}\|w_+ + w_-\|_p\\
 &= & \|G\|_{p+}
\end{array}
\]
where we used the external positivity of $G$ in deriving the first inequality.  
The above inequality clearly shows $\|G\|_p\le \|G\|_{p+}$.  
Then, since $\|G\|_p\ge \|G\|_{p+}$ trivially holds, 
we came to the conclusion $\|G\|_p= \|G\|_{p+}$.  
\end{proofof}
%

\section{Proof of \rth{th:main1}}
\label{sec:ap1}

For the proof of \rth{th:main1}, we need the next lemma.  
\begin{lemma}\cite{Ebihara_ECC2022}
Suppose $P_a\in\bbS^{n+\np}$ satisfies \rec{eq:L2+COPnew} with
$Q_a\in \COP^{\np+\nw}$.  Then, we have
\[
\scalebox{0.9}{$
  \begin{array}{@{}l}
   P_a\in \left\{P+\begin{bmatrix}
		    0_{n,n} & 0_{n,\np} \\
		    0_{\np,n} & P_p \\
		   \end{bmatrix}:\ P\in \PSD^{n+\np},\ P_p\in \COP^{\np}\right\}.  
  \end{array}$}
\]
\label{le:structure}
\end{lemma}
Now we are ready to prove \rth{th:main1}.  

\begin{proofof}{\rth{th:main1}}
For the augmented system $G_a$, we consider the trajectory of its state $x_a$
for the input $w\in L_{2+}$ with $\| w \| _2 = 1$.  
From \rec{eq:L2+COPnew}, we readily see
\[
\scalebox{0.9}{$
\begin{array}{@{}l}
\begin{bmatrix}
x_a(t) \\ w(t)
\end{bmatrix}
^\mathrm{T}
\begin{bmatrix}
P_aA_a + A_a^T P_a + C_a^T C_a & P_aB_a + C_a^T D_a \\
B_a^T P_a + D_a^T C_a & D_a^T D_a- \gamma^2 I_{\nw}
\end{bmatrix}
\begin{bmatrix}
x_a(t) \\ w(t)
\end{bmatrix}\\
+
\begin{bmatrix}
x_p(t) \\ w(t)
\end{bmatrix}^T
Q_a
\begin{bmatrix}
x_p(t) \\ w(t)
\end{bmatrix}
\le 0\ (\forall t \ge 0).  
\end{array}$}
\]
From this inequality and \rec{eq:Ga}, we have
\[
\begin{array}{@{}l}
\frac{d}{dt} \left(x_a(t)^T P_a x_a(t)\right) + z(t)^T z(t) - \gamma^2  w(t)^T w(t) \\
 +
z_p(t)^T Q_a z_p(t)
\le 0\ (\forall t \ge 0).  
\end{array}
\]
By integration over $[0,T]$, we arrive at
\begin{equation}
\begin{array}{@{}l}
\displaystyle
 x_a(T)^T P_a x_a(T) + \int_{0}^T z(t)^T z(t) - \gamma^2  w(t)^T w(t) dt\hspace*{-20mm}\\
\displaystyle
 +
\int_{0}^T z_p(t)^T Q_a z_p(t) dt
\le 0\ (\forall T > 0).  
\end{array}
\label{eq:int}
\end{equation}
Since $Q_a \in \COP^{\np+\nw}$ and $z_p\in L_{2+}$, we first note that
\[
\int_{0}^T z_p(t)^T Q_a z_p(t) dt\ge 0\ (\forall T>0).   
\]
On the other hand, 
since $x_p$ is nonnegative in $x_a(=[\ x^T\ x_p^T\ ]^T)$, 
we see from \rle{le:structure} that 
\[
x_a(T)^T P_a x_a(T) \ge 0\ (\forall T>0) .  
\]
With these facts in mind, we take the limit $T\to\infty$ in \rec{eq:int}
and obtain
\[
 \int_0^{\infty}  z(t)^T z(t) \:dt - \gamma^2 \int_0^{\infty} w(t)^T w(t) \:dt \le 0.  
\]
This clearly shows $ \| z \|_2 ^2 \le \gamma^2 \| w \|_2 ^2 = \gamma^2$.  
To summarize, we arrive at the conclusion that
\[
\| G \| _{2+} = \sup_{w \in L_{2+}, \| w \| _2 = 1} \| z \| _2 \ \le \ \gamma.  
\]
This completes the proof.  
\end{proofof}
%

\section{Proof of \rth{th:ULBA}}
\label{sec:ap2}
We start from showing the next result.  
\begin{lemma}
For the uniform infimum $\nu_p^\star\ (p\in [1,\infty),\ p=\infty)$ characterized  by \rec{eq:nups}, we have
\begin{equation}
 \nu_p^\star\ge 2^{\frac{1-p}{p}}\ (p\in[1,\infty)),\quad
 \nu_\infty^\star\ge \frac{1}{2}.  
\label{eq:lower} 
\end{equation}
\label{le:lower}
\end{lemma}
For the proof of this lemma, the next result is useful.  
\begin{lemma}
For $p\in[1,\infty)$, suppose
$x_1,x_2\in \bbR_+$ satisfies $x_1^p+x_2^p=1$.  
Then we have
\begin{equation}
x_1+x_2\le 2^{\frac{p-1}{p}}.  
\label{eq:convex}
\end{equation}
\label{le:convex}
\end{lemma}
\begin{proofof}{\rle{le:convex}}
The function $f(x)=x^p\ (p\in[1,\infty))$ is convex for $x\in\bbR_+$.  
Therefore we have
\[
 \left(\frac{x_1+x_2}{2}\right)^p \le \frac{x_1^p+x_2^p}{2}=\frac{1}{2}.  
\]
It follows that \rec{eq:convex} holds.  
\end{proofof}
\begin{proofof}{\rle{le:lower}}
For $w\in L_p$ with $\|w\|_p=1$, 
let us define $w_+,w_-\in L_{p+}$ such that $w=w_+-w_-$ by 
$w_+(t)=\max(w(t),0_n)$, 
$w_-(t)=\max(-w(t),0_n)\ (t\in[0,\infty))$.  
Then, for $p\in[1,\infty)$, we have
\begin{equation}
 \begin{array}{@{}lcl}
 \|Gw\|_p &=& \|G(w_+-w_-)\|_p\\
 &=& \|Gw_+-Gw_-\|_p\\
 &\le & \|G w_+\|_p+\|G w_-\|_p\\
 &\le & \|G\|_{p+}\|w_+\|_p+\|G\|_{p+}\|w_-\|_p\\
 &= & \|G\|_{p+}(\|w_+\|_p+\|w_-\|_p).
\end{array}
\label{eq:basic}
\end{equation}
Here, from $\|w_+\|_p^p+\|w_-\|_p^p=\|w\|_p^p=1$ and \rle{le:convex}, we have
\[
\|w_+\|_p+\|w_-\|_p\le 2^{\frac{p-1}{p}}.  
\]
From this inequality and \rec{eq:basic}, we obtain
\[
 \sup_{w\in L_p,\ \|w\|_p=1} \|Gw\|_{p}\le \|G\|_{p+} 2^{\frac{p-1}{p}}.  
\]
or equivalently, $\|G\|_{p}\le \|G\|_{p+} 2^{\frac{p-1}{p}}$.  
This clearly shows that $\nu_p^\star\ge 2^{\frac{1-p}{p}}\ (p\in[1,\infty))$.  

Similarly, for $w\in L_\infty$ with $\|w\|_\infty=1$, we have
\[
 \begin{array}{@{}lcl}
 \|Gw\|_\infty&\le & \|G w_+\|_\infty+\|G w_-\|_\infty\\
 &\le & \|G\|_{\infty+}\|w_+\|_\infty+\|G\|_{\infty+}\|w_-\|_\infty\\
 &= & \|G\|_{\infty+}(\|w_+\|_\infty+\|w_-\|_\infty)\\
 &\le & 2\|G\|_{\infty+}.  
\end{array}
\]
From this inequality, we readily obtain
\[
 \sup_{w\in L_\infty,\ \|w\|_\infty=1} \|Gw\|_{\infty}\le 2\|G\|_{\infty+}
\]
or equivalently, $\|G\|_{\infty}\le 2 \|G\|_{p+}$.   
This clearly shows $\nu_\infty^\star\ge \frac{1}{2}$.  
\end{proofof}

We next prove that \rec{eq:Gstar} holds.  To this end, the next result is useful.  
\begin{lemma}
For $w_1,w_2\in L_{p+}$, we have
\begin{equation}
\begin{array}{@{}l}
 \|w_1-w_2\|_p\le \left(\|w_1\|_p^p+ \|w_2\|_p^p\right)^{\frac{1}{p}}\ (p\in [1,\infty)),\\
 \|w_1-w_2\|_\infty\le \max(\|w_1\|_\infty,\|w_2\|_\infty).  
\end{array}
\label{eq:comps}
\end{equation}
\label{le:comps}
\end{lemma}
\begin{proofof}{\rle{le:comps}}
Let us define $\wmax\in L_{p+}$ by 
\[
\wmax(t) :=\max(w_1(t),w_2(t))\ (t\in[0,\infty)).  
\]
Then, it is clear that 
\[
\begin{array}{@{}l}
  \|w_1-w_2\|_p\le \|\wmax\|_p\ (p\in [1,\infty)),\\ 
  \|w_1-w_2\|_\infty\le \|\wmax\|_\infty.  
\end{array}
\]
On the other hand, from the definition of the $L_p$ norm of signals, 
it is clear that 
\[
\begin{array}{@{}lcl}
  \|\wmax\|_p^p& \le & \|w_1\|_p^p+\|w_2\|_p^p\ (p\in [1,\infty)),\\ 
  \|\wmax\|_\infty&\le & \max(\|w_1\|_\infty+\|w_2\|_\infty).  
\end{array}
\]
It follows that \rec{eq:comps} holds.  
\end{proofof}

We now move on to the proof \rec{eq:Gstar}.    

\begin{proofof}{\rec{eq:Gstar}}
We first prove that
\begin{equation}
 \|G^\star\|_p=2\ (p\in[1,\infty),\ p=\infty).  
\label{eq:Gstarp}
\end{equation}
To this end, note that
\begin{equation}
  \|G^\star\|_p\le \|1\|_p+\|e^{-Ls}\|_p=2\ (p\in[1,\infty),\ p=\infty).   
\label{eq:Gstarpub}
\end{equation}
To prove 
\begin{equation}
  \|G^\star\|_p\ge 2\ (p\in[1,\infty),\ p=\infty),   
\label{eq:Gstarplb}
\end{equation}
let us consider the input signal $w_N^\star\in L_p\ (p\in[1,\infty),\ p=\infty)$
defined by
\[
 w_N^\star(t)=\left\{
 \begin{array}{ll}
  0 & t< 0,\\
  1 & 2mL\le t < (2m+1)L, \\
  -1 & (2m+1)L \le t < 2(m+1)L, \\
   0 & 2NL \le t.  
 \end{array}\right.  
\]
where $N\in\bbN$ and $m=0,\cdots,N-1$.  
Then, the corresponding output $z_N^\star\in L_p\ (p\in[1,\infty),\ p=\infty)$
of the system $G^\star$ is given by
\[
 z_N^\star(t)=\left\{
 \begin{array}{ll}
  0 & t< 0,\\
  1 & 0\le t < L, \\
  2w_N^\star & L \le t < 2NL,\\
  1 & 2NL \le t < (2N+1)L,\\
  0 & (2N+1)L \le t.  
 \end{array}\right.  
\]
For the signals $w_N^\star,z_N^\star\in L_p\ (p\in[1,\infty),\ p=\infty)$, 
we readily see that
\[
\begin{array}{@{}ll}
 \|w_N^\star\|_p=(2NL)^{\frac{1}{p}}, &  \|w_N^\star\|_\infty=1,\\
 \|z_N^\star\|_p=\left((2N-1)L2^p+2L\right)^{\frac{1}{p}}, & \|z_N^\star\|_\infty=2.  
\end{array}
\]
From these results, it is obvious that $\|G^\star\|_\infty\ge 2$. 
In addition, by letting $N\to \infty$, 
we see $\|G^\star\|_p\ge 2\ (p\in[1,\infty))$.  
It follows that \rec{eq:Gstarplb} holds.  
From \rec{eq:Gstarpub} and \rec{eq:Gstarplb}, we can readily conclude that
\rec{eq:Gstarp} holds.  

We are now in a right position to to prove \rec{eq:Gstar}.  
To this end, let us consider any nonnegative input signal 
$w\in L_{p+}\ (p\in[1,\infty),\ p=\infty)$ with $\|w\|_p=1$ for $G^\star$.  
If we define $w_L\in L_{p+}\ (p\in[1,\infty),\ p=\infty)$ with $\|w_L\|_p=1$ by 
\[
 w_L(t)=
\left\{
 \begin{array}{ll}
  0 & 0 \le t < L,\\
  w(t-L) & L \le t, \\
 \end{array}
 \right.
\]
we see that the output $z$ corresponding to the input $w$ is given by $z=w-w_L$.  
By applying \rle{le:comps}, we have
\[
 \|z\|_p\le 2^{\frac{1}{p}},\ \|z\|_\infty\le 1.  
\]
Namely, for any nonnegative signal $w\in L_{p+}\ (p\in[1,\infty),\ p=\infty)$ with $\|w\|_p=1$, 
the above inequalities hold.  These results, together with, \rec{eq:Gstarp} lead us to
\begin{equation}
\begin{array}{@{}l}
\|G^\star\|_{p+}\le 2^{\frac{1}{p}}=2^{\frac{1-p}{p}}\|G\|_p=\nu_p^\star \|G\|_p\ (p\in[1,\infty)),\\
\|G^\star\|_{\infty +}\le 1=\frac{1}{2}\|G\|_\infty=\nu_\infty^\star \|G\|_\infty.  
\end{array}
\label{eq:Gstarpp}
\end{equation}
It follows from \rec{eq:lower} and \rec{eq:Gstarpp} that \rec{eq:Gstar} holds.  
\end{proofof}

We are now ready to prove \rth{th:ULBA}.  

\begin{proofof}{\rth{th:ULBA}}
Since we have verified \rec{eq:lower} and \rec{eq:Gstar},
\rth{th:ULBA} has been validated.  
\end{proofof}

\end{document}